\documentclass[12pt]{iopart}

\usepackage{amssymb}
\usepackage{graphicx}
\usepackage{fullpage}
\usepackage{psfrag}
\usepackage{color}
\usepackage{cite}
\usepackage[dvipsnames]{xcolor}

\newtheorem{lemma}{Lemma}
\newtheorem{definition}{Definition}

\newcommand{\Z}{\mathbb{Z}}
\newcommand{\R}{\mathbb{R}}
\newcommand{\N}{\mathbb{N}}
\renewcommand{\Re}{\mathrm{Re}}
\renewcommand{\Im}{\mathrm{Im}}
\newcommand{\Poincare}{Poincar\'e }
\newcommand{\out}{\mathrm{out}}
\newcommand{\iin}{\mathrm{in}}
\renewcommand{\e}{\mathrm{e}}
\newcommand{\nt}{p}
\newcommand{\toinf}{\rightarrow\infty}
\renewcommand{\vec}[1]{\mathbf{#1}}
\newcommand{\mod}{\mathrm{mod}}
\newcommand{\text}[1]{\mathrm{#1}}

\newcommand{\lmax}{\lambda_\mathrm{max}}
\newcommand{\wmax}{w^\mathrm{max}}
\newcommand{\ltwo}{\lambda_2}

\newcommand{\net}{\hat{\Sigma}}

\newcommand{\green}[1]{{{{#1}}}}

\graphicspath{{../figures/}}

\begin{document}

\title[Cycling behaviour in Rock--Paper--Scissors--Lizard--Spock]{Stability of cycling behaviour near a heteroclinic network model of Rock--Paper--Scissors--Lizard--Spock}

\author{Claire M Postlethwaite$^1$\footnote{Corresponding author (c.postlethwaite@auckland.ac.nz)} and Alastair M Rucklidge$^2$}

\address{$^1$ Department of Mathematics, University of Auckland, Private Bag 92019, Auckland 1142, New Zealand}
\address{$^2$ School of Mathematics, University of Leeds, Leeds LS2 9JT, UK}

\submitto{Nonlinearity}


\begin{abstract}The well-known game of Rock--Paper--Scissors can be used as a simple model of competition between three species. When modelled in continuous time using differential equations, the resulting system contains a heteroclinic cycle between the three equilibrium solutions representing the existence of only a single species. The game can be extended in a symmetric fashion by the addition of two further strategies (`Lizard' and `Spock'): now each strategy is dominant over two of the remaining four strategies, and is dominated by the remaining two. The differential equation model contains a set of coupled heteroclinic cycles forming a heteroclinic network. In this paper we carefully consider the dynamics near this heteroclinic network. We develop a technique to use a previously defined definition of stability (known as \emph{fragmentary asymptotic stability}) in numerical continuation software. We are able to identify regions of parameter space in which arbitrarily long periodic sequences of visits are made to the neighbourhoods of the equilibria, which form a complicated pattern in parameter space. 

\end{abstract}


\section{Introduction}

The well-known game of Rock--Paper--Scissors can be used as a simple model of competition between three species. This system has been studied extensively, in many different contexts, such as evolutionary game theory and biology~\cite{ML75,Kerr2002,Sinervo1996}, stochastic models~\cite{Kerr2002,Frey2010} and spatially dependent models~\cite{Reichenbach2007a,Frey2010,Szczesny2013, Szczesny2014}. For reviews see~\cite{Szolnoki2014,szolnoki2020pattern}. When modelled in continuous time using ordinary differential equations, the resulting system contains a heteroclinic cycle between the three equilibrium solutions representing the existence of only a single species~\cite{ML75}. This same heteroclinic cycle has also been studied in the context of fluid dynamics~\cite{BH80} and equivariant dynamical systems~\cite{dos1984, GH88}. 

The addition of two further strategies --- where each of the now five strategies  is dominant over two of the remaining four strategies, and is dominated by the remaining two --- extends the Rock--Paper--Scissors system into a network of possible states. The game of Rock--Paper--Scissors--Lizard--Spock was popularised by the TV show `The Big Bang Theory' in 2012~\cite{Cendrowski2008} (and credited to Sam Kass~\cite{Kass1995}), but equivalent networks have been around for much longer: the Wuxing cycle of five phases from Chinese philosophy has been in existence since the second or first century BCE~\cite{sivin1990science}. In modern literature, the earliest reference we can find to this network of interactions in an equivariant dynamical systems context is by Field and Richardson in 1992~\cite{Field1992}.

\begin{figure}
\begin{center}
\setlength{\unitlength}{1mm}
\begin{picture}(82,82)(0,0)
\put(0,0){\includegraphics[width=80mm]{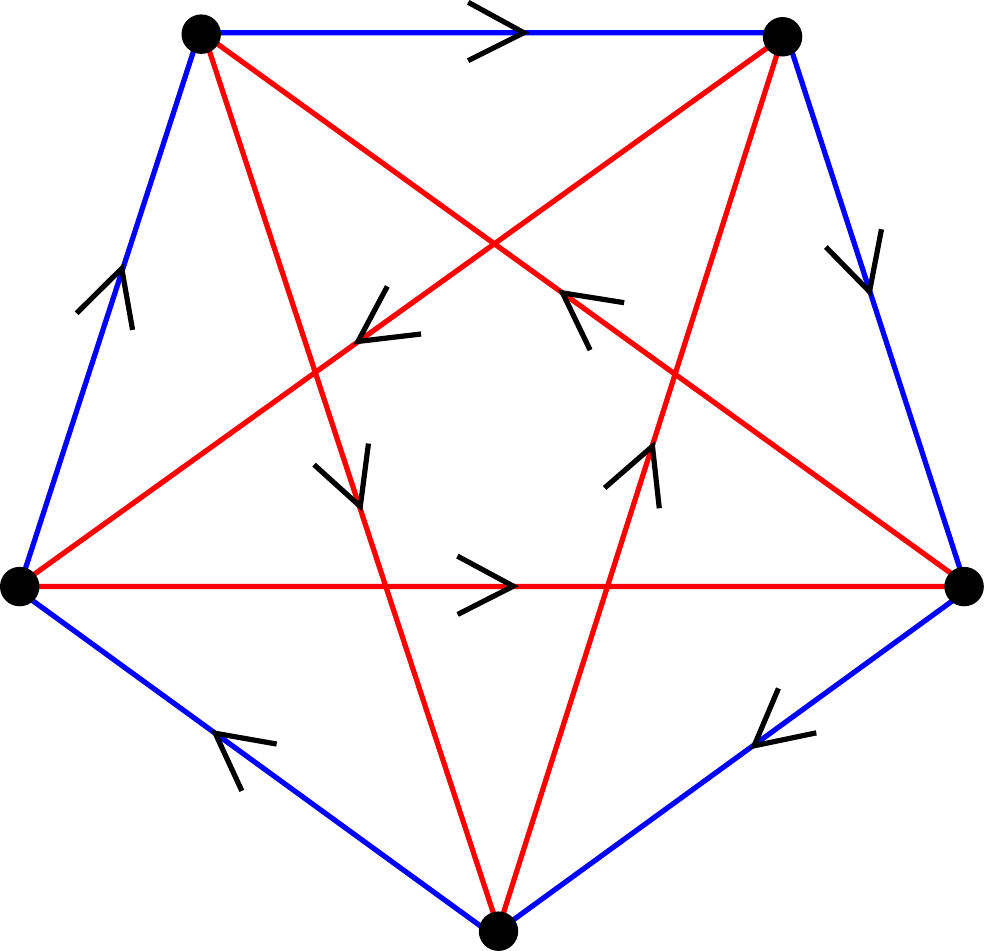}}
\put(65,77){$\xi_1$}
\put(67,70){$e_A$}
\put(62,64){$c_B$}
\put(53,71){$e_B$}
\put(56,76){$c_A$}

\put(78,33){$\xi_2$}
\put(43,0){$\xi_3$}
\put(-1,34){$\xi_4$}
\put(12,77){$\xi_5$}
\end{picture}
\end{center}
	\caption{The network of one-dimensional heteroclinic connections $\hat{\Sigma}$, between equilibria $\xi_1,\dots\xi_5$. Eigenvalues are shown near $\xi_1$, the remainder can be deduced by symmetry. The connection orbits coloured blue are of `Type A', and those coloured red are of `Type B'. Later, in section~\ref{sec:pmap} we give a formal definition of these.
	\label{fig:network}}
\end{figure}
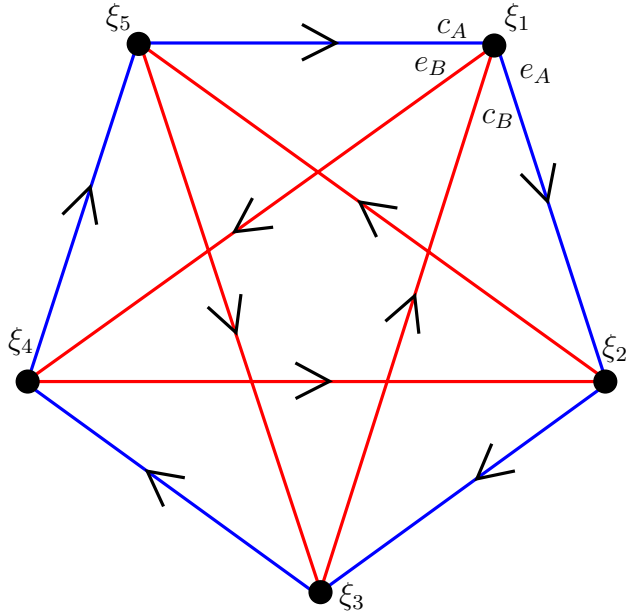

In this paper we present the first comprehensive description of the dynamics of this system when described by a system of ordinary differential equations (ODEs), namely the equations
\begin{eqnarray}
\frac{dx_{1}}{dt} & = x_1\left(1-X-c_A x_2+ e_B x_3 -c_B x_4+ e_A x_5\right) \nonumber \\
\frac{dx_{2}}{dt}& = x_2\left(1-X-c_A x_3+ e_B x_4 -c_B x_5+ e_A x_1\right) \nonumber \\
\frac{dx_{3}}{dt} & = x_3\left(1-X-c_A x_4+ e_B x_5 -c_B x_1+ e_A x_2\right) \label{eq:odes} \\
\frac{dx_{4}}{dt} & = x_4\left(1-X-c_A x_5+ e_B x_1 -c_B x_2+ e_A x_3\right) \nonumber \\
\frac{dx_{5}}{dt} & = x_5\left(1-X-c_A x_1+ e_B x_2 -c_B x_3+ e_A x_4\right) \nonumber
\end{eqnarray}
where $x_j\in\R$,  $X=x_1+x_2+x_3+x_4+x_5$, and $c_A, c_B, e_A, e_B>0$ are parameters. These ODEs contain a heteroclinic network between equilibrium solutions, shown schematically in figure~\ref{fig:network}. The stability of the entire network has been studied by Podvigina~\cite{Podvigina2020}, and Afraimovich and colleagues~\cite{Afraimovich2016}, but they do not discuss the dynamics that occur as trajectories approach the network. Some results in a stochastic setting can be found in~\cite{vukov2013, kang2013}.

The vector field generated by equations~\eref{eq:odes} (which have only linear and quadratic terms), is often thought of as being of `Lokta--Volterra' (or `May--Leonard') type. It
is topologically equivalent, in the positive orthant, to one with $\Z_2^5\rtimes \Z_5$ symmetry which can be generated by the coordinate transformation $x_j\rightarrow x_j^2$. Vector fields of this latter type, with symmetry group $\Z_2^k$ or $\Z_2^k\rtimes \Z_k$ (in both cases acting on $\R^k$) can contain heteroclinic cycles, and some have been studied in the context of equivariant bifurcation theory - for example see~\cite{GH88} for $k=3$, and~\cite{FS91} for $k=4$. In all these cases the invariant sphere theorem applies~\cite{F96}, and so all trajectories are attracted to an invariant $(k-1)$-sphere. The results we present in this paper apply only to the dynamics of~\eref{eq:odes} in the positive orthant, and as such, are independent of whether we consider these equations, or the ones with the coordinates transformed to squared variables.

For the equivariant equations (i.e.~those with linear and cubic terms), the existence and stability of heteroclinic cycles was studied in detail  for $k=4$ by Field and Swift~\cite{FS91}. The results of Field and Swift are complemented by studies of heteroclinic \emph{networks} in four dimensions by Brannath~\cite{B94} and Kirk and Silber~\cite{KS94}. Krupa and Melbourne generalised the stability results in the $k=3$ and $k=4$ cases to include cases without the $\Z_k$ permutation symmetry~\cite{KM04}. In addition to the so-called \emph{edge cycles} between equilibria with just one species present, the results of Field and Swift~\cite{FS91} show that \emph{face cycles} between equilibria with \emph{two} non-zero components can occur (see also~\cite{field2017patterns}). In the $k=5$ case, Field and Richardson~\cite{Field1992} show that these face cycles can co-exists with edge cycles, and there can also exist cycles between equilibria with \emph{three} non-zero components. In this paper, we do not consider the face-cycles, because our choice of parameters (namely, that $c_A, c_B, e_A, e_B>0$) does not allow it, but we do consider the cycle between equilibria with three non-zero components (3-face cycles in the terminology of~\cite{field2017patterns}).


Methods for determining the stability properties of an isolated robust heteroclinic cycle
between equilibria are well-established~\cite{Podvigina2012,Podvigina2013,Podvigina2017,Chossat97,KM95,KM04,Mel91,Postlethwaite2010,scheel1992},
and their implementation is generally straightforward, at least in principle, because
there is only a single route around the cycle. 
Heteroclinic cycles which are proper subsets of a heteroclinic network cannot be asymptotically stable, because there must be some points on an unstable manifold of at least one of the equilibria in the cycle which are attracted to a different part of the network. This lead to the introduction of weaker notions of stability, firstly \emph{essential asymptotic stability}, introduced by Ian Melbourne in 1991~\cite{Mel91}, and later \emph{fragmentary asymptotic stability}, introduced by Olga Podvigina in 2012~\cite{Podvigina2012}. Both of these types of stability require that the object only attract trajectories from a subset of its neighbourhood (a precise definition is given in section~\ref{sec:background} below).
 In this paper, we are not only interested in the stability of subcycles of the heteroclinic network, but also of the existence of trajectories which approach the network following a particular sequence of equilibria: this sequence may include visiting the same equilibrium multiple times but, for instance, visiting two different equilibria afterwards.

The stability results of Podvigina~\cite{Podvigina2012} are for so-called \emph{Type Z} heteroclinic cycles. The equivariant version of the vector field for~\eref{eq:odes}  is in this class. Specifically, all of the heteroclinic connections lie in fixed-point subspaces, are all of the same dimension, and the appropriate isotypic decomposition is into one-dimensional components. For more details of the group theoretic details, see definition 8 in~\cite{Podvigina2012}.

There have been some recent results giving sufficient conditions on the asymptotic stability of classes of heteroclinic network~\cite{Afraimovich2016,Podvigina2020}, which include the network we consider in this paper. However, it turns out that the network can still have very strong attracting properties even when these conditions are not satisfied. There are several other results on the stability of heteroclinic networks, for examples see~\cite{Castro2014,B94,Castro2010,Driesse2009,KS94,Kirk2010,KM95a,Postlethwaite2005,podvigina2019}, but  these are, in general, partial results and confined to specific examples. One source of difficulty is that there may be many different routes by which a trajectory can traverse a heteroclinic network, and keeping track of all possibilities in the stability calculations can be challenging. Furthermore, there must be at least one equilibrium in a network for which the unstable manifold is two-dimensional: again, keeping track of trajectories that travel close to all parts of this manifold can make the computations very involved: see~\cite{AC98,ashwin1998cycling,ashwin2004cycling,Kirk2010,Kirk2012} for examples.

The main contributions of this paper are as follows. Firstly, we give an explicit method for computing the stability, not just of sub-cycles of a heteroclinic network, but of arbitrarily complex (repeating) sequences of visits to the equilibria of the network. That is, we ask whether trajectories which visit neighbourhoods of the equilibria in the network in a particular order (the list of which can be arbitrarily long, but must eventually be repeating) are attracted to, or repelled from, the network.  \green{Using the terminology of a recent preprint from Podvigina~\cite{podvigina2021behaviour}, each of these sequences is an \emph{omnicycle}.} Secondly, we are able to adapt the conditions for heteroclinic cycle stability given by Podvigina~\cite{Podvigina2012} for use with continuation software to compute boundaries of stability mentioned above. In order to do this, we need to write the stability conditions in such a way that they generate a scalar function, continuous in the parameters of interest. Finally, we find a region of parameter space in which there is a complicated set of stability regions for increasingly complex sequences in which trajectories visit equilibria; these regions are reminiscent of the sausage-shaped resonance tongues seen in piecewise continuous systems~\cite{wei1987,campbell1996,szalai2009,simpson2016border,simpson2018structure}. A complete explanation of the shape of these regions is beyond the scope of this paper, but leaves us with many open questions and future work.

We note that there have been many other papers which prove the existence of complicated cycling behaviour of trajectories close to heteroclinic networks, and switching between different sub-cycles of the network, (e.g.~\cite{Postlethwaite2005,Kirk2010,Kirk2012,aguiar2004dynamics,homburg2010switching,castro2016switching}). To the best of our knowledge this is the first paper to give regions of parameter space where different cycling behaviours are stable (and can be observed in numerical simulations), and to explain how that stability is lost.

The remainder of this paper is organised as follows. In section~\ref{sec:background}, we give the necessary terminology and definitions required throughout the paper. In section~\ref{sec:overview} we give an overview of the dynamics of equations~\eref{eq:odes}, including the stability regions of various equilibrium solutions. We also show that in addition to the heteroclinic network between the single-population equilibria, there may also exist a heteroclinic cycle between equilibria each of which have three species present; we find the existence region and boundaries of stability for this cycle. In section~\ref{sec:pmap} we derive a \Poincare map which captures the behaviour of trajectories close to the heteroclinic network, regardless of which route between equilibria around the network is taken by a trajectory. The derivation of the \Poincare map is complicated because we have to use a mixture of Cartesian and polar coordinates in order to capture the whole of the two-dimensional unstable manifold.

\begin{figure}
\begin{center}
\setlength{\unitlength}{1mm}
\begin{picture}(130,130)(0,0)
\put(0,0){\includegraphics[trim= 0.cm 5cm 0cm 4.5cm,clip=true,width=130mm]{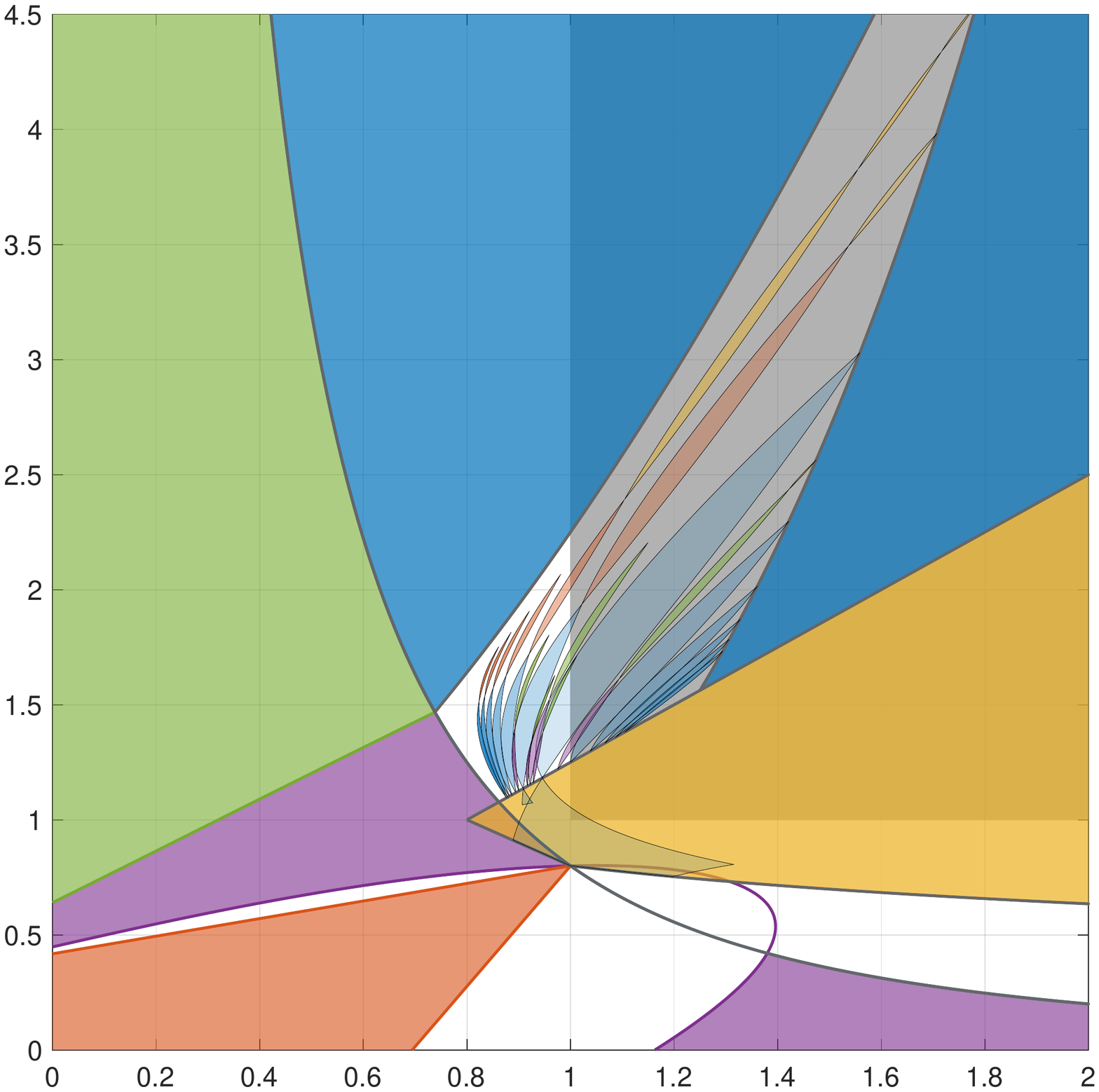}}

\put(120,0){$c_A$}
\put(2,115){$c_B$}

\put(22,73){$\xi_T$}
\put(38,30){$\Sigma_{TQ}$}
\put(93,9){$\Sigma_{TQ}$}
\put(38,12){$\xi_Q$}

\put(52,86){$AAB$}
\put(98,36){$A$}
\put(112,83){$B$}

\end{picture}
\end{center}
	\caption{Stability boundaries of various types of behaviour, in $c_A$-$c_B$ parameter space, with $e_A=1$ and $e_B=0.8$. The labels in each region correspond to stable objects: these are described in detail in later sections, but briefly, $\xi_T$ and $\xi_Q$ are equilibria, $\Sigma_{TQ}$ is a heteroclinic cycle and $A$, $B$ and $AAB$ are various patterns of approaching the heteroclinic network $\Sigma$.
	 The green, purple and orange shaded regions are the regions in which $\xi_T$, $\Sigma_{TQ}$ and $\xi_Q$ (as labelled) are asymptotically stable and the boundaries are the coloured curves: blue: $\delta_T=1$; green: $\lambda_4=0$; purple: $\delta_{TQ}=1$, orange: stability boundary of $\xi_Q$. These are described in section~\ref{sec:overview}.  The set of tongues between the regions labelled $AAB$, $A$ and $B$ are regions of fragmentary asymptotic stability of more complicated ways of approaching the network $\Sigma$. These are described in more detail in sections~\ref{sec:pmap} and~\ref{sec:num}. The grey shaded region is that for which the sufficient conditions for asymptotic stability of the network $\Sigma$ given in~\cite{Podvigina2020} and~\cite{Afraimovich2016} apply.
	\label{fig:stab_bounds}}
\end{figure}

We use 
\emph{transition matrices} in section~\ref{sec:tmatrices} to analyse the map, and apply results from Podvigina~\cite{Podvigina2012} to determine the existence of trajectories which approach the network while visiting equilibria in repeating patterns of arbitrary length. This allows us to find regions in parameter space where different patterns are stable. A summary of these results can be seen in figure~\ref{fig:stab_bounds}. The sausage-like tongues in the centre of the figure are shown in more detail in figures~\ref{fig:stab_subcycs} and~\ref{fig:stab_subcycs_zoom}, and each one represents a different pattern of visiting the equilibria in the network. The grey shaded region is that for which the sufficient conditions for asymptotic stability of the heteroclinic network given in~\cite{Podvigina2020} and~\cite{Afraimovich2016} apply. We conjecture that the grey (uncoloured) region is actually filled with infinitely many sausages for different patterns: obviously we could only compute a finite number of these. In addition, we note that there are regions of parameter space which are outside of the grey region in which the network is still strongly attracting. The numerical results and the algorithm used to generate them are discussed in detail in section~\ref{sec:num}. Finally, in section~\ref{sec:irreg} we discuss regions of parameter space in which the network appears to be attracting, but no regular pattern of visiting equilibria can be found. This irregular behaviour appears to arise in at least two different ways. Section~\ref{sec:disc} concludes.

\section{Background and definitions}
\label{sec:background}

We briefly summarise some notions and give some definitions that are used throughout this paper.
 Consider a system of ordinary differential equations
\begin{equation}\label{eq:ode1}
\dot{x}=f(x),\quad x\in\mathbb{R}^n
\end{equation}

\begin{definition}
A \emph{heteroclinic cycle} is a finite collection of equilibria $\{\xi_1, \dots, \xi_m\}$ of~\eref{eq:ode1}, together with a set of heteroclinic connections $\{\gamma_1(t),\dots, \gamma_m(t)\}$, where $\gamma_j(t)$ is a solution of~\eref{eq:ode1} such that  $\lim_{t\rightarrow -\infty} \gamma_j(t)=\xi_j$, $\lim_{t\rightarrow \infty} \gamma_j(t)=\xi_{j+1}$ and $\xi_{m+1}\equiv\xi_1$.
\end{definition}

A \emph{heteroclinic network} is a connected union of heteroclinic cycles. More generally, heteroclinic cycles and networks may connect invariant objects more complicated than equilibria, such as periodic orbits~\cite{Kirk2008} or chaotic sets~\cite{A97}, but we do not consider these possibilities here. In generic systems, heteroclinic cycles and networks are of high co-dimension, but when~\eref{eq:ode1} contains invariant subspaces, then they may exist for open sets of parameters values, that is, they are \emph{robust}.

We now give some notions of stability. The notion of fragmentary asymptotic stability was introduced by Podvigina in~\cite{Podvigina2012}. Let $X\subset\R^n$ be a set which is invariant under~\eref{eq:ode1}, and let $N_{\epsilon}$ be an $\epsilon$-neighbourhood of $X$, that is:
\[
N_{\epsilon}(X)=\{ x\in\R^n ~|~ |x-X|<\epsilon \}
\]
For $\vec{x_0}\in\R^n$, let $F_t(\vec{x_0})$ denote the flow generated by~\eref{eq:ode1}, i.e., the solution $\vec{x}(t)$ to the initial value problem starting at $\vec{x}(0) = \vec{x_0}$. We denote the $\delta$-local basin of attraction of $X$ as $\mathcal{B}_\delta(X)$:
\begin{equation}\label{eq:deltabasin}
\mathcal{B}_\delta(X)=\{x\in\R^n ~|~ | F_t(x), X|<\delta\ \forall\ t\geq 0,\ \mathrm{and}\ \lim_{t\rightarrow \infty} |F_t(x),X|=0 \}.
\end{equation}

\begin{definition}
An invariant set $X$ is \emph{asymptotically stable} if, for all $\delta>0$, there exists an $\epsilon>0$ such that
\[
N_{\epsilon}(X)\subset \mathcal{B}_\delta(X)
\]
\end{definition}

\begin{definition}[from~\cite{Podvigina2012}]
An invariant set $X$ is \emph{fragmentarily asymptotically stable (f.a.s.)} if, for any $\delta>0$, 
\[
\mu(\mathcal{B}_\delta(X))>0,
\]
where $\mu$ is the Lebesgue measure of a set in $\R^n$.
\end{definition}

If a set $X$ is fragmentarily asymptotically stable, but not asymptotically stable, then not all the points which are arbitrarily close to $X$ are attracted to $X$ in forward time. In the context of heteroclinic networks and cycles, this usually arises because a cusp- (or anticusp-) shaped region of phase space which abuts the set $X$ is excluded from the basin of attraction of $X$.

\section{Overview of the ODE system}
\label{sec:overview}

We now give an overview of the dynamics of the system of equations~\eref{eq:odes}, including the calculation of the stability of some invariants sets. Some of these are objects are subsets of the heteroclinic network which is studied for the remainder of this paper, and some are not part of this heteroclinic network.
Equations~\eref{eq:odes} are equivariant under the action of the group 
$\Gamma=\Z_5$,
generated by the element $\rho$ which has the action
 \[
 \rho(x_1,x_2,x_3,x_4,x_5)=(x_5,x_1,x_2,x_3,x_4).
 \]
 In addition, each coordinate axis, coordinate plane and three- and four-dimensional hyperplane is invariant under the flow. We label the two-, three- and four-dimensional invariant subspaces respectively  as
 \begin{eqnarray*}
 P_{jk}&=\{(x_1,x_2,x_3,x_4,x_5)| x_i=0, i\neq j, k \}, \\
 P_{jkl}&=\{(x_1,x_2,x_3,x_4,x_5)| x_i=0, i\neq j, k,l \}, \\
  P_{jklm}&=\{(x_1,x_2,x_3,x_4,x_5)| x_i=0, i\neq j, k,l ,m\}.
 \end{eqnarray*}
 
 We label the equilibria of~\eref{eq:odes} with exactly one non-zero component as $\xi_1,\dots,\xi_5$, so $\xi_1=(1,0,0,0,0)$, $\xi_2=(0,1,0,0,0)$, etc, and note that $\rho\xi_1=\xi_2$. The eigenvalues of $\xi_1$ are $-1$, $e_A$, $-c_B$, $e_B$ and $-c_A$, with eigenvectors in the $x_1$, $x_2$, $x_3$, $x_4$ and $x_5$ directions respectively (see figure~\ref{fig:network}). 
  Since we choose the parameters $c_A, c_B, e_A$ and $e_B$ to be positive, it can be shown that each of the ten two-dimensional invariant subspaces $P_{jk}$ contains either a one-dimensional heteroclinic connection from $\xi_j$ to $\xi_k$, or a one-dimensional heteroclinic connection from $\xi_k$ to $\xi_j$.  The resulting network of these one-dimensional connections is shown schematically in figure~\ref{fig:network}.

The positive orthant is invariant under the flow of~\eref{eq:odes}, and throughout this paper we consider the dynamics restricted to $\R^5_+$:
\[
\mathbb{R}^5_+=\{(x_1,x_2,x_3,x_4,x_5)~|~ x_j\geq 0, j=1,\dots, 5 \}.
\]
Within $\R^5_+$, each equilibrium $\xi_j$ has a two dimensional unstable manifold; we define the heteroclinic network $\Sigma$ to be the union of the equilibria $\xi_j$ and their unstable manifolds, namely
\[
\displaystyle \Sigma=\cup_{j=1}^5 \{ \xi_j \cup W^u(\xi^j) \}\cap \mathbb{R}^5_+
\]
We also give a name to the network of one-dimensional connections:
\[
\hat{\Sigma}=\Sigma \cap \left(\cup_{j,k=1}^5 P_{jk}\right)
\]

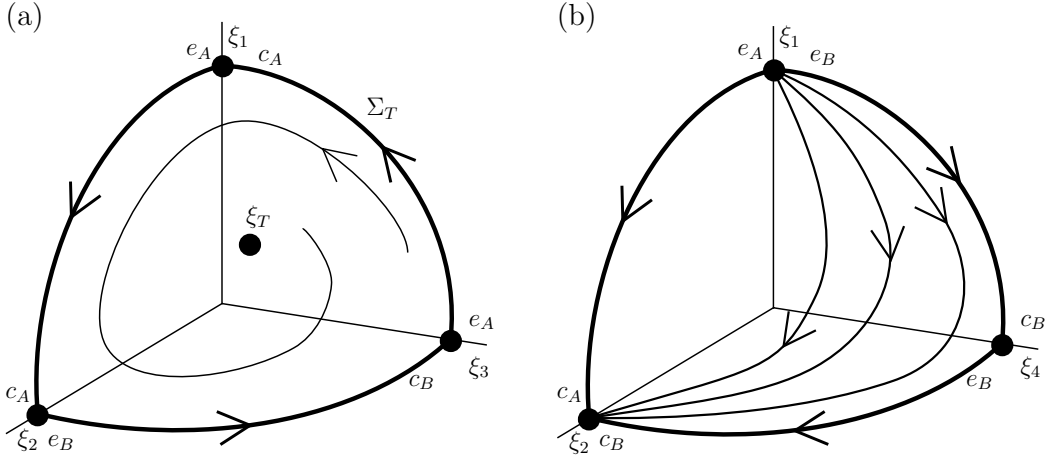
\begin{figure}
\begin{center}
\setlength{\unitlength}{1mm}

\scalebox{0.8}{
\begin{picture}(82,70)(0,0)
\put(0,0){\includegraphics[width=80mm]{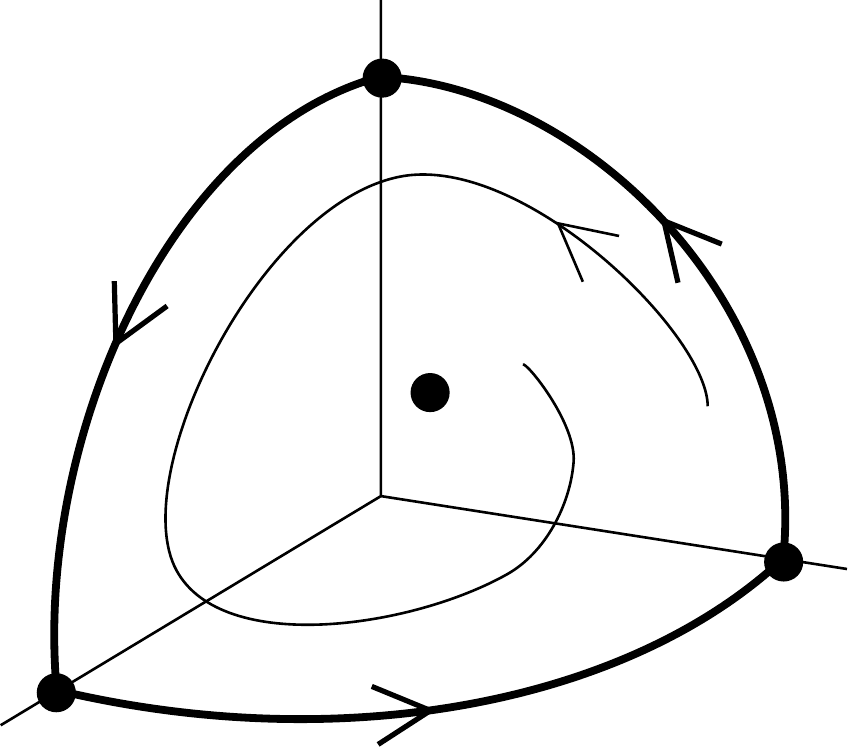}}
\put(37,67){$\xi_1$}
\put(77,12){$\xi_3$}
\put(2,0){$\xi_2$}
\put(40,37){$\xi_T$}

\put(77,20){$e_A$}
\put(67,10){$c_B$}

\put(42,64){$c_A$}
\put(30,65){$e_A$}

\put(7,0){$e_B$}
\put(0,8){$c_A$}

\put(60,55){$\Sigma_{T}$}

\put(0,70){\scalebox{1.2}{(a)}}

\end{picture}
}
\quad
\scalebox{0.8}{
\begin{picture}(82,70)(0,0)
\put(0,0){\includegraphics[width=80mm]{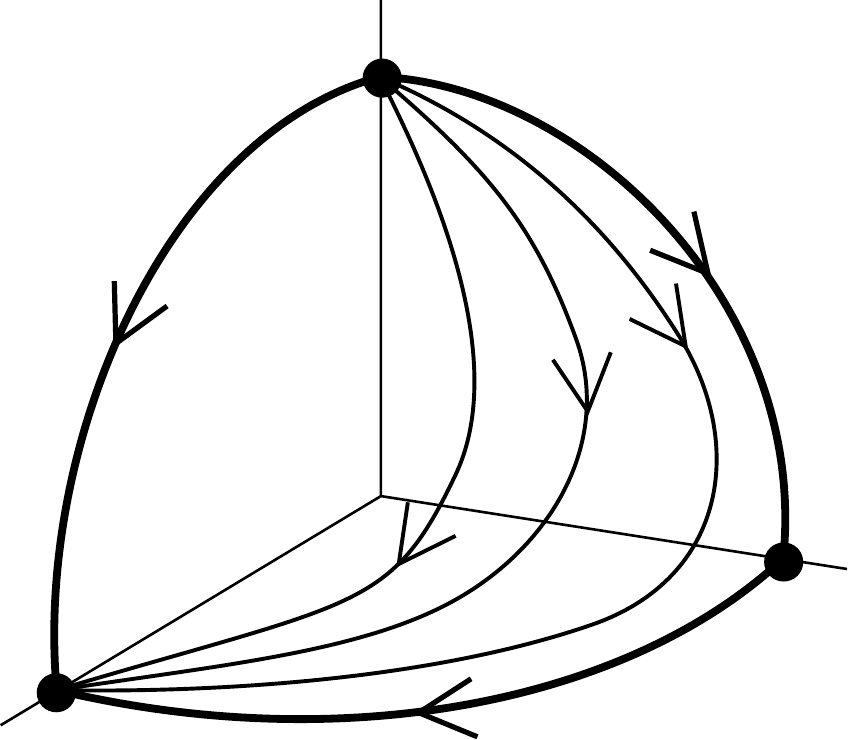}}
\put(37,67){$\xi_1$}
\put(77,12){$\xi_4$}
\put(2,-0.5){$\xi_2$}

\put(77,20){$c_B$}
\put(68,10){$e_B$}

\put(42,64){$e_B$}
\put(30,65){$e_A$}

\put(7,0){$c_B$}
\put(0,8){$c_A$}
\put(0,70){\scalebox{1.2}{(b)}}

\end{picture}
}

\end{center}
	\caption{Panel (a) shows a schematic of the dynamics of~\eref{eq:odes} restricted to the subspace $P_{123}$; the cycle $\Sigma_T$ is shown in bold. Panel (b) shows the same restricted to the subspace $P_{124}$. Eigenvalues are shown at each equilibrium $\xi_j$. The one-dimensional connections between the equilibria shown in bold (in both panels) are part of the network $\hat{\Sigma}$; the network $\Sigma$ includes the two-dimensional manifold of connections betwen $\xi_j$ and $\xi_{j+1}$.
	\label{fig:3d}}
\end{figure}

There are 10 three-dimensional invariant subspaces $P_{jkl}$ in which two of the co-ordinates are zero. The dynamics restricted to these subspaces falls into two classes, and is either symmetric (under a power of $\rho$) to the dynamics in $P_{123}$ or $P_{124}$. The dynamics of~\eref{eq:odes} restricted to each of these subspaces is shown schematically in figure~\ref{fig:3d}. Figure~\ref{fig:3d}(a) shows the dynamics restricted to $P_{123}$: it contains an equilibrium with three non-zero coordinates (labelled $\xi_T$), and a heteroclinic cycle labelled $\Sigma_T$. We discuss the dynamics in this subspace further in section~\ref{sec:P123}. Figure~\ref{fig:3d}(b) shows the dynamics restricted to $P_{124}$: notice the two-dimensional unstable manifold of $\xi_1$. 
More specifically, notice that there is a one-dimensional connection $\xi_j\rightarrow \xi_{j+3}$, but a two-dimensional connection $\xi_j\rightarrow \xi_{j+1}$ (throughout this paper, all subscripts on equilibria, subspaces, \Poincare sections and similar objects are taken $\mathrm{mod}~5$). The non-cyclic triangle of heteroclinic connections in $P_{124}$ is called a $\Delta$-clique in~\cite{ashwin2020almost}.


\begin{figure}
\begin{center}
\setlength{\unitlength}{1mm}
\begin{picture}(130,185)(-3,3)
\put(0,135){\includegraphics[trim= 3.cm 5cm 2.5cm 6.5cm,clip=true,width=130mm]{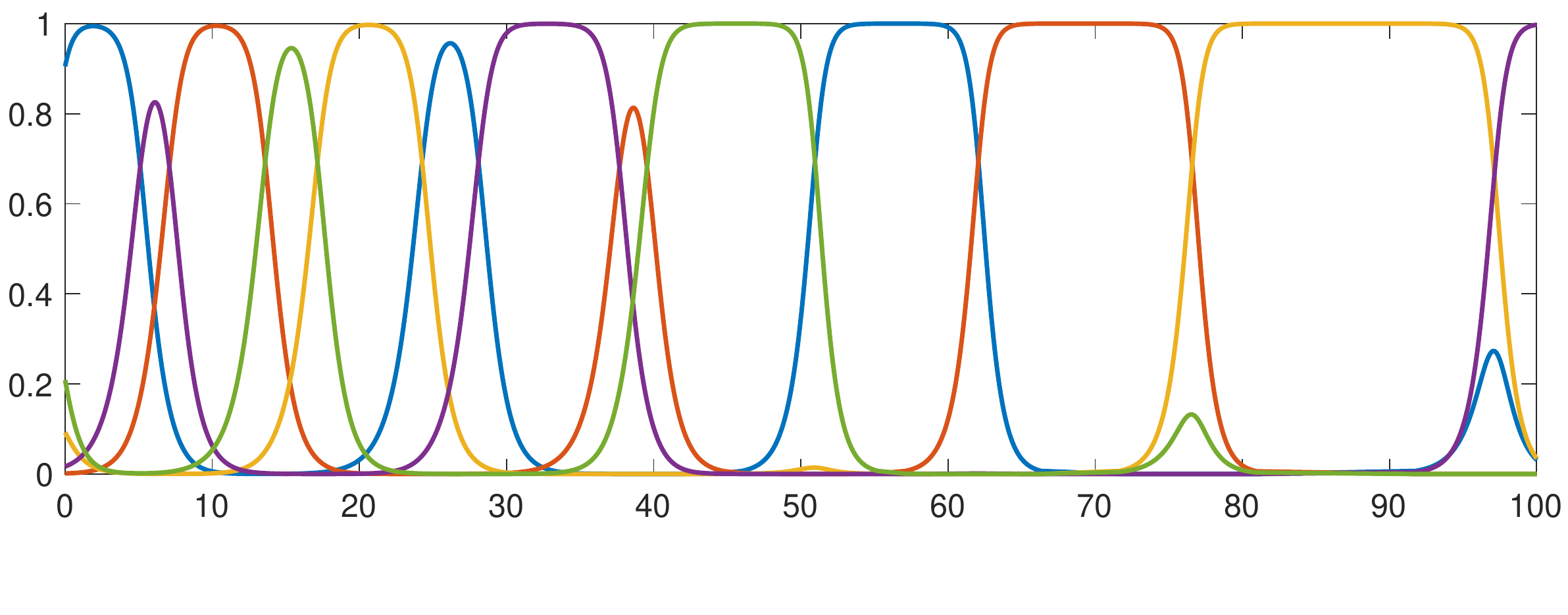}}
\put(0,90){\includegraphics[trim= 3.cm 5cm 2.5cm 6.5cm,clip=true,width=130mm]{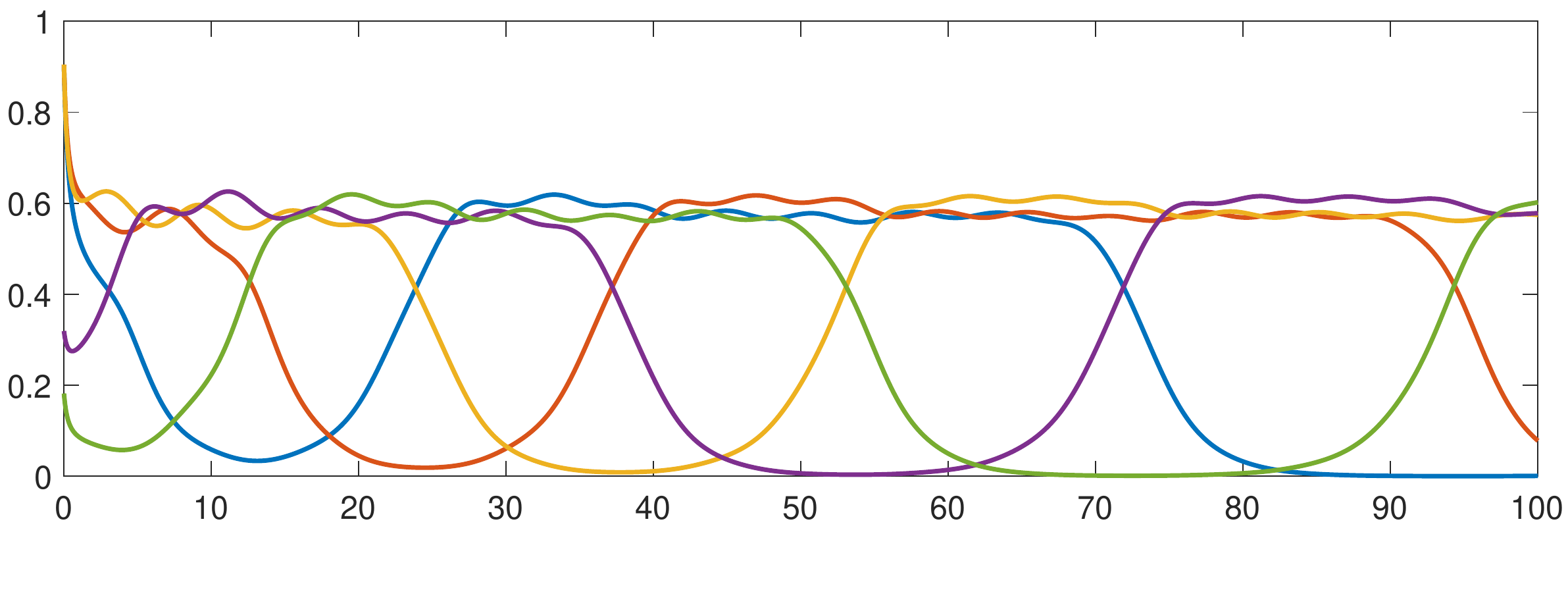}}
\put(0,45){\includegraphics[trim= 3.cm 5cm 2.5cm 6.5cm,clip=true,width=130mm]{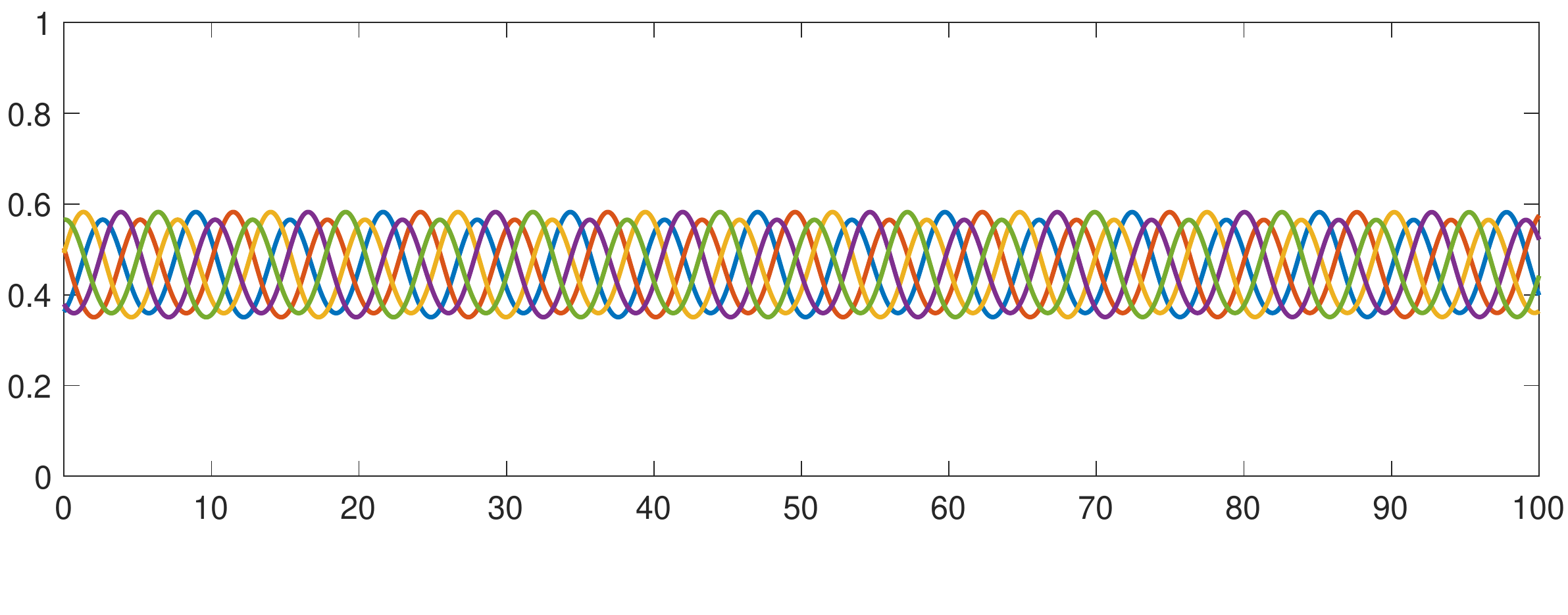}}
\put(0,0){\includegraphics[trim= 3.cm 5cm 2.5cm 6.5cm,clip=true,width=130mm]{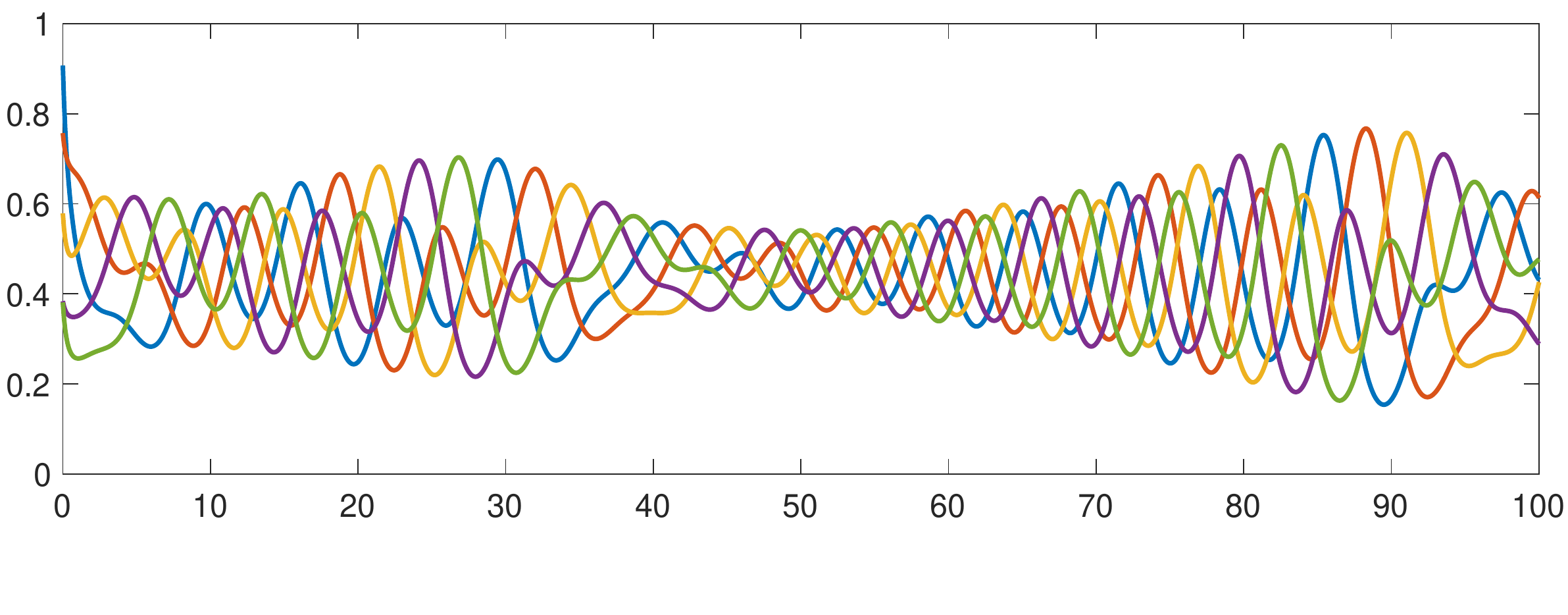}}

\put(-3,185){(a)}
\put(-3,140){(b)}
\put(-3,95){(c)}
\put(-3,50){(d)}

\put(-3,165){$x_j$}
\put(-3,120){$x_j$}
\put(-3,75){$x_j$}
\put(-3,30){$x_j$}

\put(130,4){$t$}

\end{picture}
\end{center}
\caption{The figures show typical time series of equations~\eref{eq:odes}. The lines coloured blue, red, yellow, purple and green are the coordinates $x_1,\dots,x_5$ respectively. Parameters are: panel (a): $c_A=1.2$, $c_B=1$; panel (b): $c_A=0.8$, $c_B=0.9$; panel (c): $c_A=0.89$, $c_B=0.5$; panel (d): $c_A=1.02$, $c_B=0.5$; $e_A=1$ and $e_B=0.8$ throughout. Further descriptions of the time series can be found in the text. 
	\label{fig:ts}}
\end{figure}
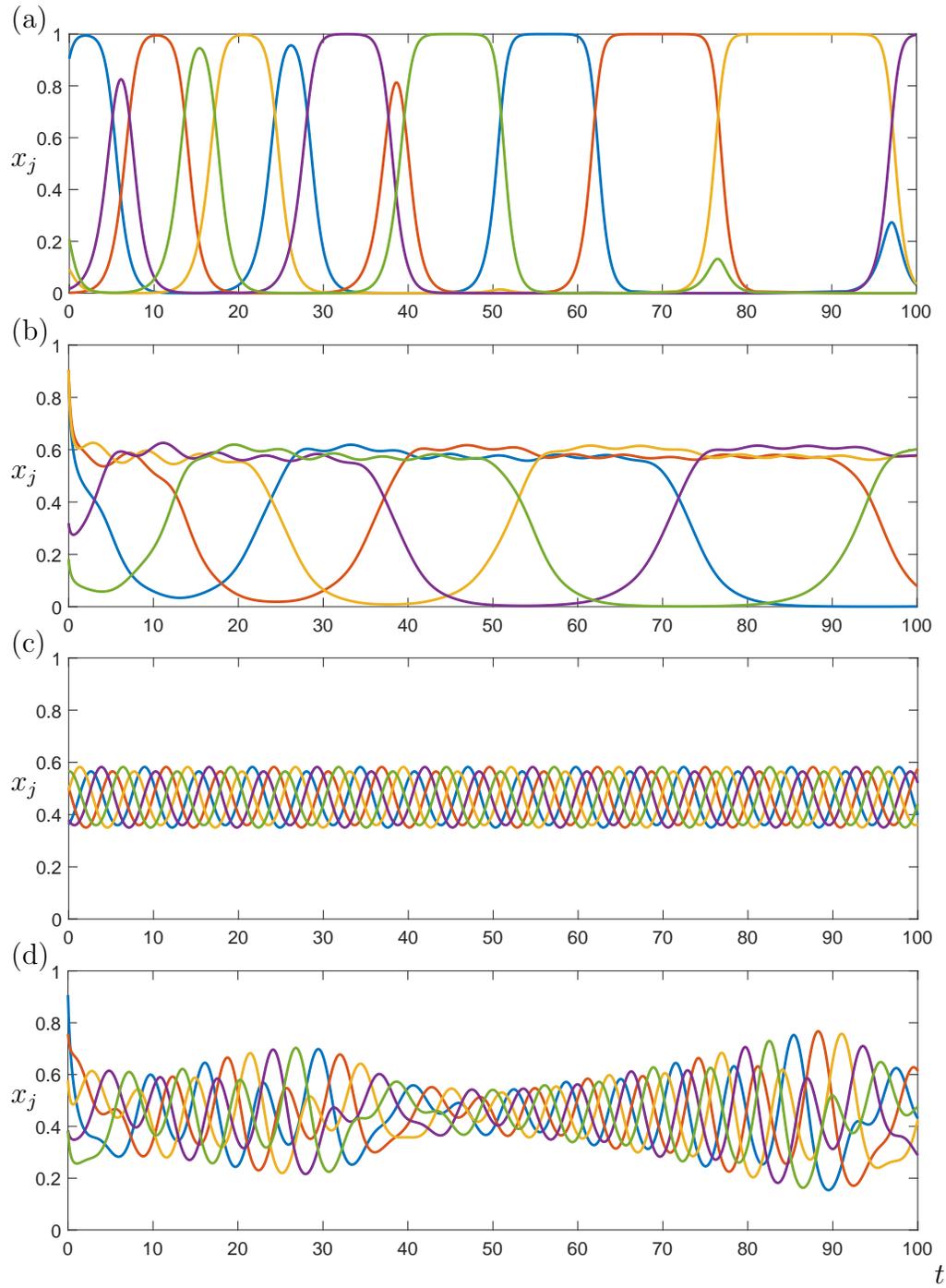

\begin{figure}
\begin{center}
\setlength{\unitlength}{1mm}
\begin{picture}(130,90)(-5,3)
\put(0,45){\includegraphics[trim= 2.8cm 5cm 2.3cm 6.5cm,clip=true,width=130mm]{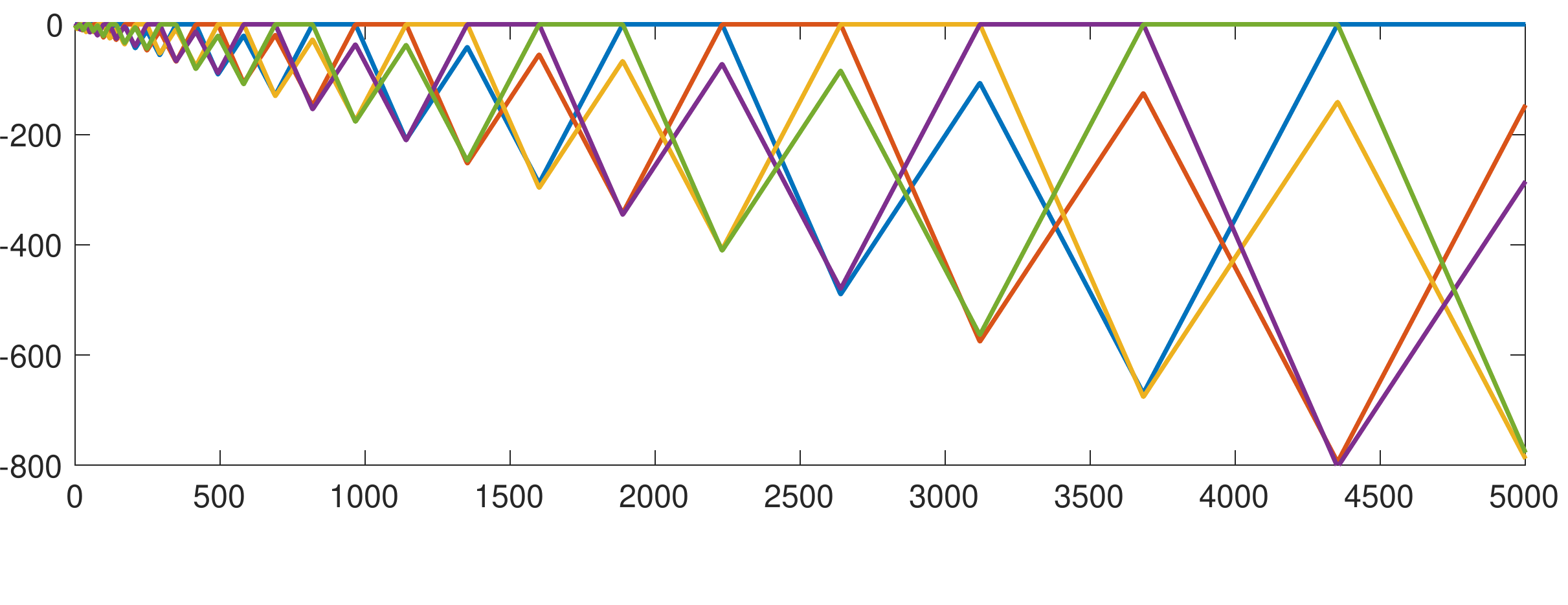}}
\put(0,0){\includegraphics[trim= 2.8cm 5cm 2.3cm 6.5cm,clip=true,width=130mm]{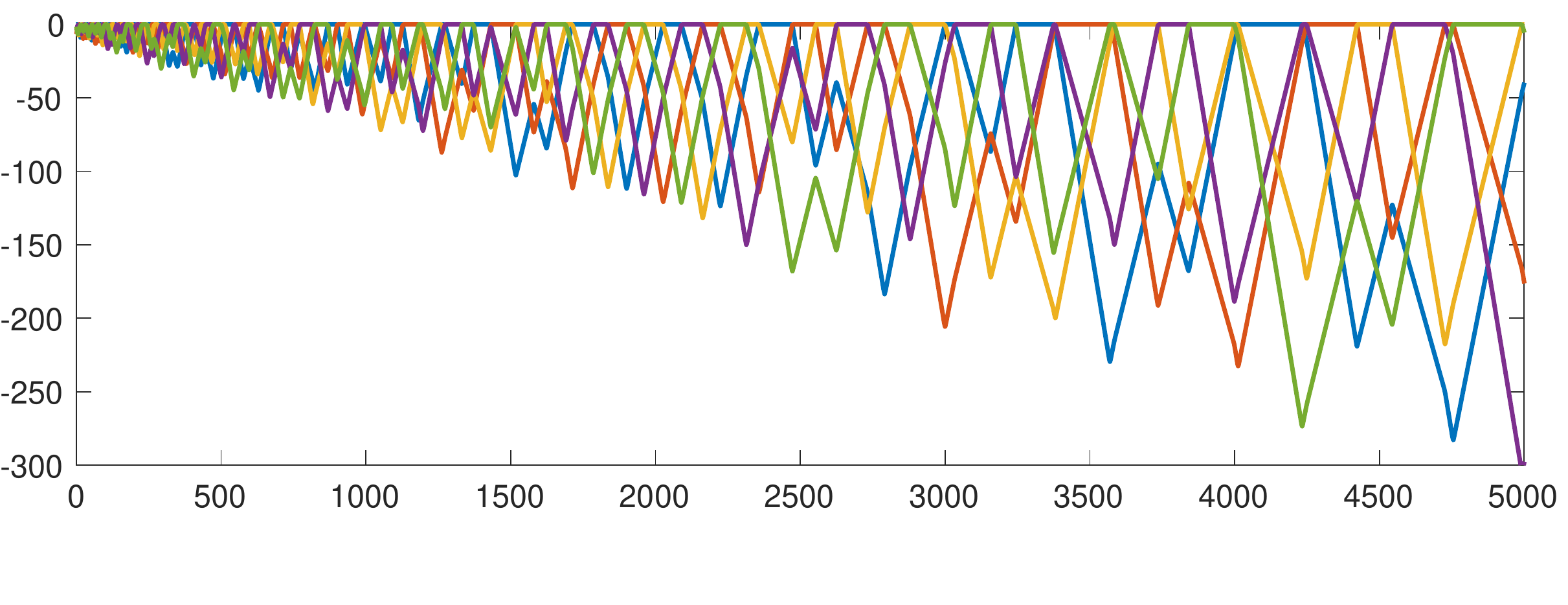}}

\put(-3,95){(a)}
\put(-3,50){(b)}
\put(-5,70){\rotatebox{90}{$\log(x_j)$}}
\put(-5,25){\rotatebox{90}{$\log(x_j)$}}
\put(130,4){$t$}

\end{picture}
\end{center}
\caption{The figures show typical time series of equations~\eref{eq:odes}, on a logarithmic scale. The lines coloured blue, red, yellow, purple and green are the logarithm of the  coordinates $x_1,\dots,x_5$ respectively. Parameters are: panel (a): $c_A=1.2$, $c_B=1$; panel (b): $c_A=1.2$, $c_B=0.7$; $e_A=1$ and $e_B=0.8$ throughout. 
	\label{fig:tslog}}
\end{figure}

In figures~\ref{fig:ts} and~\ref{fig:tslog} we show some typical time series of trajectories of equations~\eref{eq:odes}, for a variety of parameter values. Figure~\ref{fig:ts}(a) shows a trajectory approaching the heteroclinic network $\Sigma$, visiting the equilibria in the order $\xi_1$, $\xi_2$, $\xi_3$, $\xi_4$, $\xi_5$. Note that the length of time spent near each equilibrium increases as time increases. Figure~\ref{fig:ts}(b) shows a trajectory approaching the heteroclinic cycle $\Sigma_{TQ}$, in which each equilibrium has three non-zero components. Again, the time spent near each equilibrium increases.
Figure~\ref{fig:ts}(c) and (d) show periodic and quasiperiodic solutions which have bifurcated from the equilibrium $\xi_Q$.
The objects $\xi_Q$ and $\Sigma_{TQ}$ are defined later in this section. In figure~\ref{fig:tslog} we show trajectories plotted on a logarithmic scale, so that the closeness of the coordinates to zero as the equilibria are approached can be clearly seen. Figure~\ref{fig:tslog}(a) shows the same trajectory as figure~\ref{fig:ts}(a), but on a longer timescale. The geometric increase of time spent near each equilibrium solution can be clearly seen. Figure~\ref{fig:tslog}(b) shows a trajectory which is approaching the heteroclinic network $\Sigma$, but the equilibria are visited in an irregular manner.


In~\cite{Podvigina2020} and~\cite{Afraimovich2016}, it is shown that a sufficient condition for the asymptotic stability of $\Sigma$ is that $\min(c_A, c_B)>\max(e_A, e_B)$. This region is shown by the grey shading in figure~\ref{fig:stab_bounds}. As we show in the remainder of this paper, the network can have very strong stability for many regions of parameter space when this condition does not hold.

In section~\ref{sec:pmap}, we compute a \Poincare map between the equilibria $\xi_j$ which includes the whole of the two-dimensional manifold, and allows us to analyse the dynamics near $\Sigma_+$. First, for completeness, we review the remaining equilibria of~\eref{eq:odes}, their stability, and other heteroclinic cycles which can exist in this system.

\subsection{Stability of $\xi_Q$}
\label{sec:xiq}

In this section we consider the stability of the equilibrium of~\eref{eq:odes} with five non-zero components. We label this equilibrium as $\xi_Q$, and its coordinates are given by
\[
\xi_Q\equiv (x,x,x,x,x),
\]
where 
\[
x=\frac{1}{5+c_A+c_B-(e_A+e_B)},
\]
The Jacobian of~\eref{eq:odes} evaluated at $\xi_Q$ is circulant, with first row equal to
\[
J_1=[-x, -(1+c_A)x, -(1-e_B)x, -(1+c_B)x, -(1-e_A)x].
\]
The eigenvalues of circulant matrices are well known, and in this case are given by
\[
\fl
\mu_j=-x(1+(1+c_A)\omega_j+(1-e_B)\omega_j^2+(1+c_B)\omega_j^3+(1-e_A)\omega_j^4),\quad j=1,\dots,5
\]
where $\omega_j=e^{\frac{2\pi i j}{5}}$ are the fifth roots of unity. Since $\Sigma_{k=0}^5 \omega_j^k=0$, we can simplify to get
\begin{equation}\label{eq:muj}
\mu_j=-x(c_A\omega_j-e_A\bar{\omega}_j+c_B\bar{\omega}_j^2-e_B\omega_j^2),\quad j=1,\dots,5 .
\end{equation}


We thus find
\begin{equation}
\Re(\mu_1)=-x\left((c_A-e_A)\cos\left(\frac{2\pi}{5}\right)-(c_B-e_B)\cos\left(\frac{\pi}{5}\right)\right),
\end{equation}
and 
\begin{equation}
\Re(\mu_2)=-x\left(-(c_A-e_A)\cos\left(\frac{\pi}{5}\right)+(c_B-e_B)\cos\left(\frac{2\pi}{5}\right)\right).
\end{equation}


Recall that $x>0$, and $\frac{\cos\left(\frac{\pi}{5}\right)}{\cos\left(\frac{2\pi}{5}\right)}=\frac{\sqrt{5}+1}{\sqrt{5}-1}=\beta>1$.
We then have that $\Re(\mu_1)<0$ if 
\[
(e_A-c_A)<\beta(e_B-c_B)
\]
and $\Re(\mu_2)<0$ if 
\[
(e_B-c_B)<\beta(e_A-c_A).
\]

For $\xi_Q$ to be stable, we need to satisfy both of these conditions, which requires $\frac{e_A-c_A}{e_B-c_B}<\beta$,and $\frac{e_B-c_B}{e_A-c_A}<\beta$. This stability boundary is shown in parameter space in figure~\ref{fig:stab_bounds} by the orange lines, and the region in which $\xi_Q$ is stable is shaded orange.
Both boundaries are Hopf bifurcations. A periodic solution resulting from one of these bifurcations is shown in figure~\ref{fig:ts}(c). Figure~\ref{fig:ts}(d) shows a quasiperiodic solution to which this solution evolves after further bifurcations.

\subsection{The subspace $P_{123}$ }
\label{sec:P123}

In this section we consider the dynamics within the subspace $P_{123}$, as shown schematically in figure~\ref{fig:3d}(a).
This subspace contains a heteroclinic cycle between $\xi_1$, $\xi_2$ and $\xi_3$. We refer to this cycle as $\Sigma_T$. There are four further symmetric copies of this sub-cycle contained in $\hat{\Sigma}$, related by (powers of) the symmetry $\rho$.
The subspace $P_{123}$ also contains an equilibrium with $x_1, x_2$ and $x_3$ non-zero, which we label as $\xi_T$ (shown also in figure~\ref{fig:3d}(a)). The full system~\eref{eq:odes} contains a further  four symmetric copies of this equilibrium, also with three non-zero components, related to $\xi_T$ again by powers of the symmetry $\rho$.

The heteroclinic cycle $\Sigma_T$ is equivalent to the frequently-studied Guckenheimer--Holmes cycle~\cite{ML75,GH88}, with the removal of the rotational symmetry relating the three equilibria. Regardless of this, the stability of $\Sigma_T$ is simple to compute using a \Poincare map, and stability results are similar. Namely, for the dynamics of~\eref{eq:odes} restricted to $P_{123}$, the equilibrium $\xi_T$ and the cycle $\Sigma_T$ exchange stability at a degenerate Hopf bifurcation when $\delta_T=\frac{c_A^2 c_B}{e_A^2 e_B}=1$. The cycle $\Sigma_T$ is asymptotically stable for the dynamics restricted to $P_{123}$ if $\delta_T>1$ and the equilibrium $\xi_T$ is stable if $\delta_T<1$. (Note that this degenerate Hopf bifurcation can be broken with the addition of higher order terms, resulting in a branch of periodic orbits.) 
 
For the full dynamics in $\R^5$, $\Sigma_T$ cannot be asymptotically stable (since, for instance, any points arbitrarily close to $\xi_1$ that lie on the heteroclinic connection between $\xi_1$ and $\xi_4$ will asymptote to $\xi_4$), but it can still be fragmentarily asymptotically stable. We consider the stability of $\Sigma_T$ in detail in section~\ref{sec:stabsigmaT}, once we have computed the \Poincare map for the network.

The equilibrium $\xi_T$ can lose stability to perturbations transverse to the subspace $P_{123}$. 
The Jacobian of $\xi_T$ is block diagonal, and so the eigenvalues of $\xi_T$ can be easily computed. We label the eigenvalues with eigenvectors in the $x_4$ and $x_5$ directions as $\lambda_4$ and $\lambda_5$, respectively. When $\delta_T<1$ and both $\lambda_4<0$ and $\lambda_5<0$, $\xi_T$ is asymptotically stable. In figure~\ref{fig:stab_bounds}, we show the curve $\delta_T=1$ in blue, and the curve $\lambda_4=0$ in green ($\lambda_5<0$ for all parameters shown in figure~\ref{fig:stab_bounds}). The region shaded green shows where $\xi_T$ is stable.

Now suppose that $\delta_T<1$, $\lambda_5<0$, but $\lambda_4>0$, and consider the  dynamics within the four-dimensional subspace $P_{1234}$. This subspace contains the equilibria $\xi_T$ and $\rho\xi_T$, and by the symmetry $\rho$, within this subspace, $\rho\xi_T$ will be a sink. $\xi_T$ will be a saddle with a one-dimensional unstable manifold (the unstable eigenvector of $\xi_T$ points in the $x_4$ direction), and so there will exist a heteroclinic connection between $\xi_T$ and $\rho\xi_T$. This connection will be robust to perturbations preserving the invariant subspaces in the system. 
Again, by the  symmetry $\rho$, there are also robust connections between $\rho \xi_T$ and $\rho^2 \xi_T$, $\rho^2\xi_T$ and $\rho^3\xi_T$, etc, creating a heteroclinic cycle between the five equilibria $\rho^j \xi_T$, $j=0,\dots ,4$; we label this heteroclinic cycle $\Sigma_{TQ}$.

 $\Sigma_{TQ}$ is not part of a heteroclinic network: each equilibrium has a one-dimensional unstable manifold, and so we can compute its stability using \Poincare maps in the usual manner (see~\cite{Postlethwaite2006} for a very similar calculation). We find that it becomes unstable as the quantity $\delta_{TQ}=-\lambda_4/\lambda_5$ decreases through $1$, shown by a purple curve in figure~\ref{fig:stab_bounds}. The stable region of $\Sigma_{TQ}$ is shaded purple in figure~\ref{fig:stab_bounds}.

\section{Dynamics near the heteroclinic network}
\label{sec:pmap}

In this section, we construct a \Poincare map which captures the behaviour of trajectories close to the heteroclinic network $\hat{\Sigma}$ (shown in figure~\ref{fig:network}) between the equilibria $\xi_1$, $\xi_2$, $\xi_3$, $\xi_4$ and $\xi_5$. We then use this map to compute the stability of any given periodic sequence of transitions around the heteroclinic network, using the results of Podvigina~\cite{Podvigina2012}.

We first define some terminology.  We are interested in trajectories of~\eref{eq:odes} that remain close to $\Sigma$, and so we can describe the motion in terms of the itinerary around the network, that is, the sequence of visits to the equilibria $\xi_j$. The following notation follows that used in~\cite{Ashwin2016}. For fixed $0<H<\frac{1}{2}\min_{i,j}|\xi_i-\xi_j|$, and $\vec{x}\in\R^5_+$, we define
$$
M(\vec{x}) := \left\{\begin{array}{rl}
k & \mbox{ if there exists a $k$ such that } |\vec{x}-\xi_k|\leq H\\
0 & \mbox{ otherwise.}\end{array}\right.
$$
Note that the choice of maximum allowed $H$ means that $M(\vec{x})$ is uniquely defined. When $M(F_t(\vec{x_0})))=k$ we say $F_t(\vec{x_0})$ is {\em close to $\xi_k$}, and note that the set of points $\vec{x}$ for which $M(\vec{x})=j$ forms a ball around $\xi_j$. 
For a trajectory $F_t(\vec{x_0})$ of~\eref{eq:odes} we define
\begin{equation}
\tilde{M}(t) = \{ M(F_{\tilde{t}}(\vec{x_0})) ~|~ \tilde{t} = \sup\{ \tilde{t}\leq t~|~ M(F_{\tilde{t}}(\vec{x_0}))\neq 0\} \}.
\label{eq:lastvisited}
\end{equation}
If $F_{t}(\vec{x_0})$ is not close to any equilibria, that is, $M(F_{t}(\vec{x_0}))=0$, then the $\tilde{t}$ in the above expression will be the time when the trajectory was \emph{most recently} close to an equilibrium, and hence $\tilde{M}(t)$
 gives the `last visited equilibrium'. If a trajectory starts close to an equilibrium at $t=0$ this will always be non-zero. The trajectory $F_t(\vec{x_0})$ can be thus characterised as an \emph{itinerary} $m(n)$ of {\em epochs} $\tau(n)$:
$$
\{ (m(n),\tau(n)) ~|~ n\in\N\}
$$
such that $\tilde{M}(t)=m(n)$ for the interval $t\in[\tau(n),\tau(n+1))$, and $m(n+1)\neq m(n)$. We say that a trajectory \emph{transitions} from $\xi_j$ to $\xi_k$ at time $t$ if there exists an $n\in\mathbb{N}$ such that $m(n)=j$, $m(n+1)=k$ and $\tau(n+1)=t$.

For the dynamics of~\eref{eq:odes}, notice that if a trajectory stays close to the heteroclinic network $\net$, it can make only two types of transitions between equilibria. The first type, which we call ``Type A'' is a transition from an equilibrium $\xi_j$ to the equilibrium $\xi_{j+1}$. 
The second type of transition is from an equilibrium $\xi_j$ to the equilibrium $\xi_{j+3}$. We label this type of transition as ``Type B''. Given any sequence of transitions between equilibria, we can thus translate the sequence $m(n)$ into a word in the alphabet $\{A, B\}$. Specifically, for a trajectory $F_t(\vec{x_0})$ with itinerary $m(n)$, we write
\[
W(F_t(\vec{x_0}))=\{w(n) ~|~ n\in\N\}
\]
where
\[
w(n)=\left\{\begin{array}{l} A \quad \mathrm{if} \quad m(n+1)-m(n) = 1~(\mathrm{mod}~5), \\ 
B  \quad \text{if} \quad m(n+1)-m(n) = 3~(\mod~5). \end{array}  \right.
\]
We say $w(n)$ is \emph{eventually periodic} with period $p$ and \emph{root sequence} $w^\star\in \{A,B\}^p$, if there exists an $N\in\mathbb{N}$ such that for all $n>N$, $w(n+p)=w(n)$, and $w(N+j)=w^{\star}(j)$ for $j=1,\dots,p$. The minimal period $\hat{p}$ of an eventually periodic sequence $w(n)$ is the smallest such $p$.
The root sequence of an eventually periodic sequence is unique up to a cyclic permutation of letters.

We give a couple of examples. If a trajectory approaches the heteroclinic cycle in the subspace $P_{123}$ between the equilibria $\xi_1$, $\xi_2$ and $\xi_3$ (called $\Sigma_T$ in section~\ref{sec:P123}), then the resulting itinerary will be eventually periodic with root sequence $AAB$. If a trajectory approaches the heteroclinic cycle between all five equilibria $\xi_1$, $\xi_2$, $\xi_3$, $\xi_4$ and $\xi_5$, in that order, the resulting itinerary will be eventually periodic with root sequence $A$. If a trajectory approaches the same five equilibria but in the order $\xi_1$, $\xi_4$, $\xi_2$, $\xi_5$, $\xi_3$, then the resulting itinerary will be eventually periodic with root sequence $B$.

\begin{definition}
The \emph{$\delta$-local basin of attraction} $\mathcal{B}_\delta (w)$ of a sequence $w\in \{A,B\}^p$, $p\in\N$, is the set of points $\vec{x}_0\in\R^5_+$ such that $W(F_t(\vec{x}_0))$ is eventually periodic with root sequence $w$, $|F_t(x), \Sigma|<\delta$ for all  $\ t\geq 0$, and $\lim_{t\rightarrow \infty} |F_t(\vec{x}_0)-\Sigma|=0$.
\end{definition}

\begin{definition}
A sequence $w\in \{A,B\}^p$, $p\in\N$, is \emph{fragmentarily asymptotically stable (f.a.s.)} if for all $\delta>0$, $\mu( \mathcal{B}_\delta (w))>0$, where $\mu$ is the Lebesgue measure.
\end{definition}

In the calculations which follow, we determine for which parameter values particular sequences are f.a.s.
When we talk about the \emph{stable region} for a particular root sequence, we mean this region of parameter space. Note that if there exists at least one root sequence which is f.a.s., then the network $\Sigma$ is also f.a.s.

Figure~\ref{fig:stab_subcycs} shows a more detailed summary of our results than figure~\ref{fig:stab_bounds}, including a complicated region of strings of sausage-shaped stability regions for different sequences (shown in more detail in figure~\ref{fig:stab_subcycs_zoom}). Recall that the region of parameter space that the sufficient condition for stability of $\Sigma$ given by Podvigina~\cite{Podvigina2020} and Afraimovich~\cite{Afraimovich2016} is $\min(c_A,c_B)>\max(e_A,e_B)$, which is $\min(c_A,c_B)>1$ for our choice of parameters $e_A$ and $e_B$: this is the region shaded gray in figure~\ref{fig:stab_bounds}. There are clearly large regions of parameters space outside this region in which at least one root sequence is attracting, which means that the network is fragmentarily asymptotically stable.

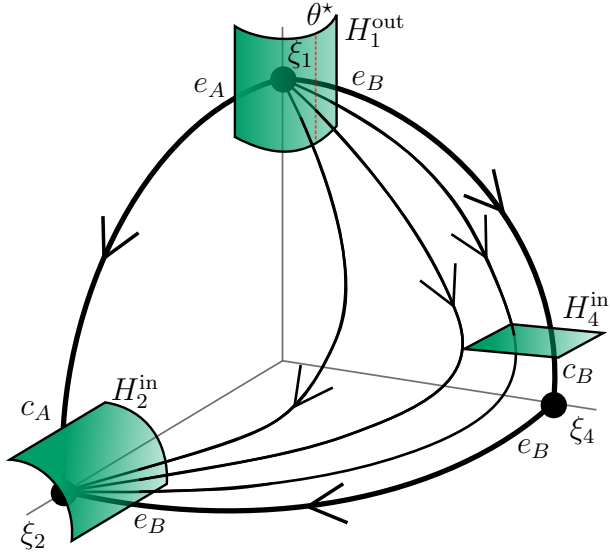
\begin{figure}
\begin{center}
\setlength{\unitlength}{1mm}
\begin{picture}(82,82)(0,0)
\put(0,0){\includegraphics[width=80mm]{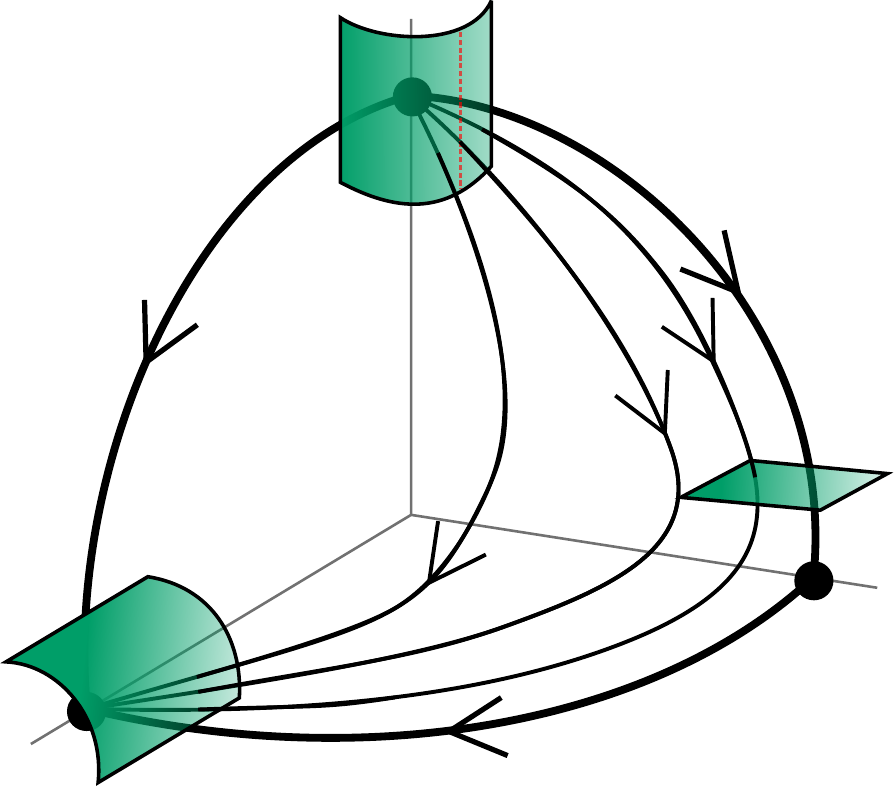}}
\put(37.5,64){$\xi_1$}
\put(75,14){$\xi_4$}
\put(2,0){$\xi_2$}

\put(74,22){$c_B$}
\put(68,12){$e_B$}

\put(46,61){$e_B$}
\put(25,60){$e_A$}

\put(17,2){$e_B$}
\put(2,17){$c_A$}

\put(45,67){$H_1^{\out}$}
\put(40,69){$\theta^{\star}$}
\put(74,30){$H_4^{\iin}$}
\put(14,19){$H_2^{\iin}$}

\end{picture}
\end{center}
	\caption{A schematic of the dynamics within the subspace $P_{124}$ and slices of the \Poincare sections $H_1^{\out}$, $H_4^{\iin}$ and $H_2^{\iin}$. The dashed line on $H_1^{\out}$ indicates the location of $\theta_{24}=\theta^{\star}$.
	\label{fig:Psections}}
\end{figure}

\subsection{Construction of a \Poincare map}

We follow a standard procedure described, for instance, in~\cite{KS94,Kirk2012}. As in~\cite{Kirk2012}, we have to use a combination of Cartesian and polar coordinates in order to capture the dynamics near the whole of the two-dimensional unstable manifold.

We begin our construction of a \Poincare map by defining \Poincare sections near $\xi_1$, for some $h\ll 1$, as follows:
\begin{eqnarray*}
H_1^{\iin}&=\{(x_1,x_2,x_3,x_4,x_5) ~|~ x_3^2+x_5^2 =h, |x_j|<h, j=2,\dots,5, |x_1-1|<h \} \\
H_1^{\out}&=\{(x_1,x_2,x_3,x_4,x_5) ~|~ x_2^2+x_4^2 =h, |x_j|<h, j=2,\dots,5, |x_1-1|<h  \}.
\end{eqnarray*}
As is usual in these types of calculations, there is an attracting invariant sphere~\cite{F96}, and thus we know that the radial direction will not play a part in the stability calculations. Discounting the radial direction then, the \Poincare sections are three-dimensional. We label points on each section as follows:
\[
(x_2, x_4, \theta_{53})\in H_1^{\iin}, \quad (x_3, x_5, \theta_{24})\in H_1^{\out}
\]
where $\tan\theta_{53}=\frac{x_5}{x_3}$ and $\tan\theta_{24}=\frac{x_2}{x_4}$. \Poincare sections near the other equilibrium are defined by the symmetry $\rho$, that is, $H_{j+1}^{\iin}=\rho H_j^{\iin}$, and $H_{j+1}^{\out}=\rho H_j^{\out}$ (subscripts taken $\mod~5$), and we show slices of $H_1^{\out}$, $H_2^{\iin}$ and $H_4^{\iin}$ schematically in figure~\ref{fig:Psections}.

We construct a local map $\phi$ from $H_1^{\iin}$ to $H_1^{\out}$ by solving the linearised flow of~\eref{eq:odes} near $\xi_1$, and assuming that the time the trajectory takes to travel from $H_1^{\iin}$ to $H_1^{\out}$ is given by $T$.

\begin{eqnarray*}
&\phi:H_1^{\iin} \rightarrow H_1^{\out} \\
&\phi(x_2^{\iin}, x_4^{\iin}, \theta_{53}^{\iin})=(x_3^{\out},x_5^{\out},\theta_{24}^{\out})
\end{eqnarray*}
where
\begin{eqnarray*}
x_3^{\out}& = h\cos \theta_{53}^{\iin} \e^{-c_B T} \\
x_5^{\out}& = h\sin \theta_{53}^{\iin} \e^{-c_A T} \\
\tan \theta_{24}^{\out} & = \frac{x_2^{\iin}}{x_4^{\iin}}  \e^{(e_A-e_B) T} \\
\end{eqnarray*}
and $T$ can be found by solving
\[
(x_2^{\iin})^2 \e^{2e_A T}+(x_4^{\iin})^2 \e^{2e_B T}=h^2
\]

The computation of the global map between equilibria is more subtle: depending on the outgoing coordinates on $H_1^{\out}$, the trajectory may hit $H_2^{\iin}$ or $H_4^{\iin}$ first. From figure~\ref{fig:Psections}, we can see that if $\theta_{24}$ is very close to $\pi/2$, we expect the trajectory to hit $H_2^{\iin}$ first, and if $\theta_{24}$ is very close to $0$, we expect the trajectory to hit $H_4^{\iin}$ first. 

More precisely, within the subspace $P_{124}$, the invariant sphere theorem means that the dynamics are restricted to a two-dimensional manifold, and so there will be some $\theta^{\star}$, such that if $\theta_{24}<\theta^{\star}$ the trajectory first hits $H_4^{\iin}$, and if $\theta_{24}>\theta^{\star}$ the trajectory first hits $H_2^{\iin}$. As we will see in the calculations that follow, the invariant subspaces force the (lowest order) global map to be diagonal in the coordinates $x_3$ and $x_5$, which means that the dynamics slightly outside of the subspace $P_{124}$ will be similar.

We thus define two global maps
\begin{eqnarray*}
\psi_{12}:& H_1^{\out}\rightarrow H_2^{\iin} \\
\psi_{14}:&  H_1^{\out}\rightarrow H_4^{\iin}
\end{eqnarray*}
We note that the relevant coordinates on $H_2^{\iin}$ are $x_3$, $x_5$ and $\theta_{14}$, (where $\tan\theta_{14}=\frac{x_1}{x_4}$) and the relevant coordinates on $H_4^{\iin}$ are $x_5$, $x_2$ and $\theta_{31}$ (where $\tan\theta_{31}=\frac{x_3}{x_1}$). 

First consider $\psi_{12}$, and write $\psi_{12}(x_3^{\out},x_5^{\out}, \theta_{24}^{\out})=(x_3^{\iin},x_5^{\iin}, \theta_{14}^{\iin})$, where $\theta_{24}^{\out}\in(\theta^{\star},\pi/2)$. Due to the invariance of the subspaces $P_{1245}$ and $P_{1234}$, it is clear that to lowest order, we will have
\begin{eqnarray*}
x_3^{\iin}&=A_3 x_3^{\out} \\
x_5^{\iin}&=A_5 x_5^{\out} 
\end{eqnarray*}
where $A_3$ and $A_5$ are order 1 `global' constants. Note also that $\theta_{24}=\pi/2$ defines an invariant subspace --- which is equivalent to $\theta_{14}=\pi/2$ (both are the subspace $P_{1235}$). Thus, if we write
\[
\theta_{14}^{\iin}=g_{12}(x_3^{\out},x_5^{\out}, \theta_{24}^{\out})
\]
then we know that $g_{12}(x_3^{\out},x_5^{\out},\pi/2)=\pi/2$, for any $x_3^{\out}$, $x_5^{\out}$. If $|\theta_{24}^{\out}-\pi/2|$ is small, we can Taylor expand $g_{12}$, and get to lowest order
\[
\frac{\pi}{2}-\theta_{14}^{\iin}=A_\theta \left(\frac{\pi}{2}-\theta_{24}^{\out}\right)
\]
where $A_\theta$ is another order 1 global constant.

Next consider the global map $\psi_{14}$, and write  $\psi_{14}(x_3^{\out},x_5^{\out}, \theta_{24}^{\out})=(x_2^{\iin},x_5^{\iin}, \theta_{13}^{\iin})$, and recall that $\theta_{24}^{\out}\in(0,\theta^{\star})$. Again, by the invariance of the subspace with $x_5=0$, we can write
\begin{eqnarray*}
x_5^{\iin}&=B_5 x_5^{\out} 
\end{eqnarray*}
for an order 1 global constant $B_5$. Since $x_2^{\iin}$ is small on $H_4^{\iin}$ (by the definition), then we can use the invariance of the subspace $x_2=0$, and Taylor expand about zero to get
\begin{eqnarray*}
x_2^{\iin}&=B_2\theta_{24}^{\out}
\end{eqnarray*}
Similarly, we find that to lowest order,
\begin{eqnarray*}
\theta_{13}^{\iin}&=B_3x_3^{\out}.
\end{eqnarray*}

\subsection{Equivariant coordinates}

We can now write down a \Poincare map from a single \Poincare section to itself, using the symmetry $\rho$ which can map $H^{\iin}_j$ to $H^{\iin}_{j+1}$. To do this, we write $H^{\iin}\equiv H^{\iin}_1$, and introduce equivariant coordinates (see, e.g.~\cite{KM04}) $x^e_A$, $x^e_B$ and $\theta_c$ on $H^{\iin}$ and, for consistency, we use $x^c_A$, $x^c_B$ and $\theta_e$ on $H^{\out}_1$. When the trajectory is close to $\xi_1$, we have the equivalencies $x^e_A\equiv x_2$, $x^e_B\equiv x_4$, $\theta_c\equiv \theta_{53}$, and $x^c_A\equiv x_5$, $x^c_B\equiv x_3$, $\theta_e\equiv \theta_{24}$.

The local map can be written in these coordinates as:
\begin{eqnarray*}
&\phi:H^{\iin} \rightarrow H^{\out}_1 \\
&\phi (x_A^{e}, x_B^e, \theta_c)=(x_A^c, x_B^c, \theta_e)
\end{eqnarray*}
where
\begin{eqnarray*}
x_B^c& = h\cos \theta_{c} \e^{-c_B T} \\
x_A^c& = h\sin \theta_{c} \e^{-c_A T} \\
\tan \theta_{e} & = \frac{x_{A}^e}{x_{B}^e}  \e^{(e_A-e_B) T} \\
\end{eqnarray*}
and $T$ can be found by solving
\begin{equation}
\label{eq:T}
(x_{A}^e)^2 \e^{2e_A T}+(x_{B}^e)^2 \e^{2e_B T}=h^2
\end{equation}

 To write the global map in these same coordinates, we can apply the symmetry $\rho$ to the \Poincare sections $H^{\iin}_2$ and $H^{\iin}_4$, namely $\rho^4 H^{\iin}_2=H^{\iin}_1$ and  $\rho^2 H^{\iin}_4=H^{\iin}_1$. In the equivariant coordinates, we then find that in a neighbourhood of $\xi_2$, we have  $x^e_A\equiv x_3$, $x^e_B\equiv x_5$, $\theta_c\equiv \theta_{14}$, and in a neighbourhood of $\xi_4$, we have  $x^e_A\equiv x_5$, $x^e_B\equiv x_2$, $\theta_c\equiv \theta_{31}$.

This gives the following equivariant global map
\begin{eqnarray*}
&\psi:H^{\out}_1\rightarrow H^{\iin} \\
&\psi(x_A^c, x_B^c, \theta_e)= (x_A^e, x_B^e, \theta_c)
\end{eqnarray*}
where, if $0<\theta_e<\theta^{\star}$, to lowest order we get:
\begin{eqnarray*}
x_A^e& = B_5 x_A^c \\
x_B^e& = B_2 \theta_e \\
\theta_c&=B_3 x_B^c
\end{eqnarray*}
and if $\theta^{\star}<\theta_e<\pi/2$, then
\begin{eqnarray*}
x_A^e& = A_3 x_A^c \\
x_B^e& = A_5 x_A^c\\
\theta_c &=g_{12}(x_B^c,x_A^c, \theta_e)
\end{eqnarray*}
where we have limited information about the function $g_{12}$. Specifically, we only know that if $|\theta_e -\pi/2|$ is small, to lowest order
\[
\theta_c=\frac{\pi}{2}- A_\theta \left(\frac{\pi}{2}-\theta_e \right).
\]
In the global map, $B_2$, $B_3$, $B_5$, $A_3$, $A_5$ and $A_{\theta}$
are order 1 global constants.

\subsection{Approximations to the \Poincare map}

In this section we consider the relevant approximations that can be made to the \Poincare map in the case where the trajectories remain always close to at least one two-dimensional invariant subspace. That is, the angles $\theta_e$ and $\theta_c$ are always either close to $0$ or close to $\pi/2$ in a neighbourhood of any of the equilibria. There are four cases in total, which depend on the coordinates on $H^{\iin}$. 

The coordinates on $H^{\iin}$ give information both  about what type of transition between equilibria will be made once the trajectory leaves $H^{\out}$, and also about which type of transition was made by the trajectory as it approached $H^{\iin}$. Namely, the relative size of $x_A^e$ and $x_B^e$ give you information about where the trajectory will go when it leaves $H^{\out}$, and whether $\theta_c$ is close to $0$ or $\pi/2$ tells you from which direction the section $H^{\iin}$ was approached.

More precisely, we consider the relative sizes of the coordinates $x_A$ and $x_B$ on $H^{\iin}$, in order to approximate the time $T$. From equation~\eref{eq:T}, we can see that if $x_A^e  \ll  {x_B^e}^\frac{e_A}{e_B}$, then to lowest order, $T=-\frac{1}{e_B} \log(x_B^e)$. In this case, the trajectory leaving $H^{\out}$ will make a Type $B$ transition (as defined in at the start of section~\ref{sec:pmap}).
We can then compose the local and global map to get a full return map on $H^{\iin}$ that we call $\Phi_B$:
\begin{equation}\label{eq:phib}
\Phi_B \pmatrix{ x_A \cr x_B \cr \theta }  = \pmatrix{   C_A \sin\theta  x_B^{\frac{c_A}{e_B}} \cr C_B x_A x_B^{-\frac{e_A}{e_B}} \cr C_\theta \cos\theta x_B^{\frac{c_B}{e_B}} } 
\end{equation}
where we have dropped the superscripts $e$ and $c$ on coordinates for clarity. 

Similarly, if $x_A ^e \gg  {x_B^e}^\frac{e_A}{e_B}$, then to lowest order, $T=-\frac{1}{e_A} \log(x_A^e)$. In this case, the trajectory leaving $H^{\out}$ will make a  Type $A$ transition. When we compose the local and global maps in this case, we find the second full return map on $H^{\iin}$, which we call $\Phi_A$: 
\begin{equation}\label{eq:phia}
 \Phi_A \pmatrix{ x_A \cr x_B \cr \theta }  = \pmatrix{ D_A \cos\theta x_A^{\frac{c_B}{e_A}}  \cr D_B  \sin\theta x_A^{\frac{c_A}{e_A}} \cr \frac{\pi}{2}- D_\theta x_Bx_A^{-\frac{e_B}{e_A}} } 
\end{equation}
(again, dropping the superscripts on the coordinates).

We now consider the composition of these maps as different routes around the network are taken. Recall that $x_A$ and $x_B$ are small. If a Type $B$ transition is made, then the resulting $\theta$ will be very small. If a Type $A$ transition is made, the resulting $\theta$ will be close to $\pi/2$. We thus define $\varphi=\pi/2-\theta$, and are able to give four different \Poincare maps, which depend on both the \emph{previous} equilibrium visited, and the one which will be visited next, that is, on both the previous and current transition types.

The final results give four complete \Poincare maps, $H^{\iin}\rightarrow H^{\iin}$ depending on the sequence of connections. Here, $\Phi_{X\rightarrow Y}$ ($X,Y\in\{A,B\}$) is the appropriate map to use when the \emph{previous} transition is of type $X$, and the \emph{current} transition is of type $Y$.
\begin{eqnarray*}
\Phi_{A\rightarrow A}:H^{\iin}\rightarrow H^{\iin}  & \quad \quad  \Phi_{A \rightarrow B}:H^{\iin}\rightarrow H^{\iin} \\
\Phi_{A\rightarrow A} \pmatrix{ x_A \cr x_B \cr \varphi }  = \pmatrix{ D_A \varphi x_A^{\frac{c_B}{e_A}}  \cr D_B  x_A^{\frac{c_A}{e_A}} \cr D_\theta x_Bx_A^{-\frac{e_B}{e_A}} } & \quad \quad
\Phi_{A \rightarrow B} \pmatrix{ x_A \cr x_B \cr \theta }  = \pmatrix{ C_A x_B^{\frac{c_A}{e_B}}  \cr C_B x_A x_B^{-\frac{e_A}{e_B}} \cr C_\theta \varphi x_B^{\frac{c_B}{e_B}} } 
  \\
  & \\
   \Phi_{B\rightarrow B}:H^{\iin}\rightarrow H^{\iin}  & \quad \quad \Phi_{B \rightarrow A}:H^{\iin}\rightarrow H^{\iin}\\
 \Phi_{B\rightarrow B} \pmatrix{ x_A \cr x_B \cr \theta }  = \pmatrix{ C_A \theta x_B^{\frac{c_A}{e_B}}  \cr C_B  x_A x_B^{-\frac{e_A}{e_B}} \cr C_\theta x_B^{\frac{c_B}{e_B}} } & \quad \quad
 \Phi_{B \rightarrow A} \pmatrix{ x_A \cr x_B \cr \varphi }  = \pmatrix{ D_A x_A^{\frac{c_B}{e_A}}  \cr D_B \theta x_A^{\frac{c_A}{e_A}} \cr D_\theta x_B x_A^{-\frac{e_B}{e_A}} } 
  \\
\end{eqnarray*}

Let the sequence of letters $Z=Z_1Z_2\dots Z_m$, where each $Z_i\in\{A,B\}$, be a root sequence of transitions for an eventually periodic sequence corresponding to a trajectory close to the network. The map given by the composition
\[
\Phi=  \Phi_{Z_{m-1} \rightarrow Z_m} \circ   \Phi_{Z_{m-2} \rightarrow Z_{m-1}} \circ \dots \circ \Phi_{Z_2 \rightarrow Z_3}  \circ \Phi_{Z_1 \rightarrow Z_2} \circ \Phi_{Z_m \rightarrow Z_1}
\]
is the return map to $H^{\iin}$ which describs the dynamics after the sequence $Z$ of transitions has been made. For any $Z$, this map has a fixed point at zero, and the stability of the zero fixed point under iterations of the map corresponds to the stability of the heteroclinic network, for trajectories following this particular sequence of transitions.
 The following section discusses how to compute the stability of the zero fixed point in maps of this type.

\subsection{Transition matrices and stability analysis}
\label{sec:tmatrices}

In order to analyse the maps given in the previous section, and apply the results of Podvigina~\cite{Podvigina2012}, we first need to define the \emph{transition matrix} of a map~\cite{KM04,FS91}. Let $G$ be the set of mappings $\Psi:\R^\nt\rightarrow\R^\nt$ that have at
lowest order the form
\[\Psi(x_1,\dots,x_\nt)=(C_1x_1^{\alpha_{11}}x_2^{\alpha_{12}}\cdots
x_\nt^{\alpha_{1\nt}},\dots, C_\nt x_1^{\alpha_{\nt1}}\cdots x_\nt^{\alpha_{\nt\nt}})\]
for real constants $\alpha_{ij}\geq 0$ and $C_i$ non-zero, $1\leq i,j \leq p$. $G$ is clearly closed under
composition. We define the \emph{transition matrix}
of $\Psi$ to be the $\nt\times \nt$ real
matrix $M(\Psi)$ with entries $[M(\Psi)]_{ij}=\alpha_{ij}$. It is easily
verified that if $\Psi_1, \Psi_2 \in G$, then
\[M(\Psi_2\circ \Psi_1)=M(\Psi_2)M(\Psi_1).\]
Any $\Psi\in G$ has a fixed point at $x_1=\dots=x_\nt=0$.
The zero fixed point of the map $\Psi$ will be stable if all the row sums of
$M(\Psi)^N$ diverge to $+\infty$ as $N\toinf$. Conversely, if any of the
row sums of $M(\Psi)^N$ tends to $0$, then the fixed point is unstable.

Furthermore, if the vector $\mathbf{v}$ (with components $v_j$) is an eigenvector of $M(\Psi)$, then the curve in $\mathbb{R}^p$ with $x_j=q^{v_j}$, $q\in [0,\infty)$, $j=1,\dots, p$, is invariant under the map $\Psi$. If the eigenvalue $\lambda$ corresponding to $\mathbf{v}$ is real and greater than 1, then along this curve, points contract towards the origin. If, furthermore, $\lambda$ is also the eigenvalue of $M(\Psi)$ with largest magnitude, then this curve is attracting under iteration of $\Psi$. Podvigina~\cite{Podvigina2012} extends these results to show that if these conditions on the eigenvalue and eigenvector of $M(\Psi)$ are satisfied, then the origin in $\Psi$ is fragmentarily asymptotically stable. That is, the basin of attraction of the origin can be extended from this one-dimensional curve to a set with positive measure. The stability of the origin in $\Psi$ is then related to the stability of the heteroclinic cycle for which the return map $\Psi$ was derived. 
A summary of the results from~\cite{Podvigina2012} is given in definition~\ref{defn:fas} and lemma~\ref{lem:stab}  below.

We now write down the transition matrices associated with the maps $\Phi_{A\rightarrow A}$, $\Phi_{A\rightarrow B}$, $\Phi_{B\rightarrow B}$ and $\Phi_{B\rightarrow A}$:
\begin{eqnarray*}
M_{A \rightarrow A}  =  \pmatrix{ \frac{c_B}{e_A} & 0 & 1 \cr
\frac{c_A}{e_A} & 0 & 0 \cr
-\frac{e_B}{e_A} & 1 & 0 }
   \quad \quad   & M_{A \rightarrow B}  =
   \pmatrix{ 0 & \frac{c_A}{e_B} & 0 \cr
1 & -\frac{e_A}{e_B} & 0 \cr
0 & \frac{c_B}{e_B} & 1 } \\
  & \\
M_{B \rightarrow B}  =  \pmatrix{ 
0 & \frac{c_A}{e_B} & 1 \cr
1 & -\frac{e_A}{e_B}  & 0 \cr
0 & \frac{c_B}{e_B} & 0 }
   \quad \quad   & M_{B \rightarrow A}  =
   \pmatrix{ \frac{c_B}{e_A} & 0 & 0 \cr
\frac{c_A}{e_A} & 0 & 1 \cr
-\frac{e_B}{e_A} & 1 & 0 }
\end{eqnarray*}

Again, let $Z=Z_1Z_2\dots Z_m$,  $Z_i\in\{A,B\}$ be a root sequence. We write $\mathcal{M}_{Z^1}$ to be the following product of $m$ transitions matrices:
\[
\mathcal{M}_{Z^1}= M_{Z_{m-1} \rightarrow Z_m} M_{Z_{m-2} \rightarrow Z_{m-1}}  \dots M_{Z_2 \rightarrow Z_3}  M_{Z_1 \rightarrow Z_2} M_{Z_m \rightarrow Z_1}
\]
and similarly also write $\mathcal{M}_{Z^j}$, for $j=2,\dots, m$ to be the product of $m$ transition matrices:
\[
\fl
\mathcal{M}_{Z^j}= M_{Z_{j-2} \rightarrow Z_{j-1}} M_{Z_{j-3} \rightarrow Z_{j-2}} \dots   M_{Z_1 \rightarrow Z_2}  M_{Z_m \rightarrow Z_1} M_{Z_{m-1} \rightarrow Z_m} \dots  M_{Z_j \rightarrow Z_{j+1}}  M_{Z_{j-1} \rightarrow Z_j}
\]
We say that the set of $m$ matrices $\{\mathcal{M}_{Z^j}, j=1,\dots m\}$ is a collection of transitions matrices describing the root sequence $Z$. When we compute the stability of a root sequence, we must consider the \Poincare maps associated with each of these $m$ transitions matrices: this is equivalent to starting the \Poincare map near each of the $m$ different equilibria in the sequence. The stability conditions (given below) depend on the eigenvalues and eigenvectors of these matrices. All of the matrices in the collection $\{\mathcal{M}_{Z^j}\}$ will have the same eigenvalues, but the eigenvectors will in general be different.

The following definition uses the conditions for stability in Lemma 5 of Podvigina's 2012 paper~\cite{Podvigina2012}. 
\begin{definition}\label{defn:fas}
Let $\mathcal{M}$ be a transition matrix for a map $g$. Let $\lambda_\mathrm{max}$ be the eigenvalue with largest absolute value of the matrix $\mathcal{M}$, and $w^{\mathrm{max}}$ be the associated eigenvector. Suppose $\lambda_\mathrm{max}\neq 1$. Then $\mathcal{M}$ is of \emph{fragmentary asymptotic stability type} if the following conditions hold:
\begin{enumerate}
\item $\lambda_\mathrm{max}$ is real
\item $\lambda_\mathrm{max}>1$
\item $w^{\mathrm{max}}_l w^{\mathrm{max}}_q>0$ for all $l,q$.
\end{enumerate}
\end{definition}
Note that the last condition is equivalent to requiring all the entries of the eigenvector $w^{\mathrm{max}}$ to be non-zero and of the same sign.

Using this definition, we then have the following result, which follows directly from Lemma 5 in Podvigina~\cite{Podvigina2012}.
\begin{lemma}
Let $\{\mathcal{M}_{Z^j}, j=1,\dots m\}$ be a collection of transition matrices for a root sequence $Z$ of length $m$. If, for each $j$, the matrix $\mathcal{M}_j$ is of fragmentary asymptotic stability type, then the root sequence $Z$ is fragmentarily asymptotically stable.
\label{lem:stab}
\end{lemma}

A sequence which is f.a.s.~can lose this stability by violating any of the three conditions listed in definition~\ref{defn:fas} for any of the transition matrices in its collection. As noted earlier, all the transition matrices in a collection will have the same eigenvalues, so if either condition \emph{(i)} or \emph{(ii)} are violated, this happens for all the transition matrices in the collection at the same parameter value. The violation of condition $1.$ corresponds to a resonance bifurcation, and would be expected to be associated with the appearance of a long-period periodic orbit. The violation of conditions \emph{(ii)} or \emph{(iii)} are not associated with the bifurcation of any other invariant objects (as noted by Podvigina~\cite{Podvigina2012}, and seen also in a particular example in~\cite{Postlethwaite2010}.) If condition \emph{(iii)} is violated, generically it will only be violated for a single one of the matrices in a collection at one time. This gives information on how instability manifests, that is, at which point in the sequence a nearby trajectory exits. We give specific examples of this when discussing the losses of stability we observe in section~\ref{sec:num}.

\subsection{Examples of stability conditions}

\subsubsection{Stability of five-cycles}

The network contains two cycles of length 5. As discussed earlier, these cycles correspond to root sequences $A$, and $B$ respectively. 
 The stability of these is computed in~\cite{Podvigina2013}, and the conditions are as follows.
The root sequence $A$ is f.a.s.~if
\[
c_A+c_B>e_A+e_B, \quad c_Ae_A>c_Be_B, \quad \text{and} \quad c_Ac_B^3>e_Ae_B^3.
\]
When $c_A+c_B=e_A+e_B$, then the largest eigenvalue of $M_{A\rightarrow A}$ becomes equal to one. This is equivalent to a resonance-type bifurcation (see, e.g.~\cite{scheel1992,Postlethwaite2010a}). One might expect the appearance of a long-period periodic orbit to be associated with this bifurcation, but we have not been able to find one numerically. We suspect that this is due to the lack of higher order terms in equations~\eref{eq:odes} (as is the case with the resonance bifurcation of $\Sigma_T$, discussed in section~\ref{sec:P123}).

 When $c_Ae_A=c_Be_B$, the third component of the eigenvector $w^{\mathrm{max}}$ is equal to zero. This means when $c_Ae_A<c_Be_B$, there is no longer an open subset of $H^{\iin}$ which returns to $H^{\iin}$ after a single iteration of the map $\Phi_{A\rightarrow A}$. When  $c_Ac_B^3=e_Ae_B^3$, a pair of complex conjugate eigenvalues of $M_{A\rightarrow A}$ has the same amplitude as the real eigenvalue with largest amplitude. Again, this means that when $c_Ac_B^3<e_Ae_B^3$,  there is no longer an open subset of $H^{\iin}$ which returns to $H^{\iin}$ after a single iteration of the map $\Phi_{A\rightarrow A}$. Both of these mechanisms of stability loss do not result in the appearance of a bifurcating object.

The root sequence $B$ is f.a.s.~if
\[
c_A+c_B>e_A+e_B, \quad c_Be_B>c_Ae_A, \quad \text{and} \quad c_A^3e_B>c_Be_A^3.
\]
The $B$ sequence loses stability in the same way for the first and the third of these conditions, respectively, as the $A$ sequence does, explained above. When $c_Be_B=c_Ae_A$, the matrix $M_{B\rightarrow B}$ has two real eigenvalues, with the same absolute value, one of which is positive and one of which is negative.

\begin{figure}
\begin{center}
\setlength{\unitlength}{1mm}
\begin{picture}(130,130)(0,0)
\put(0,0){\includegraphics[trim= 0.cm 5cm 0cm 4.5cm,clip=true,width=130mm]{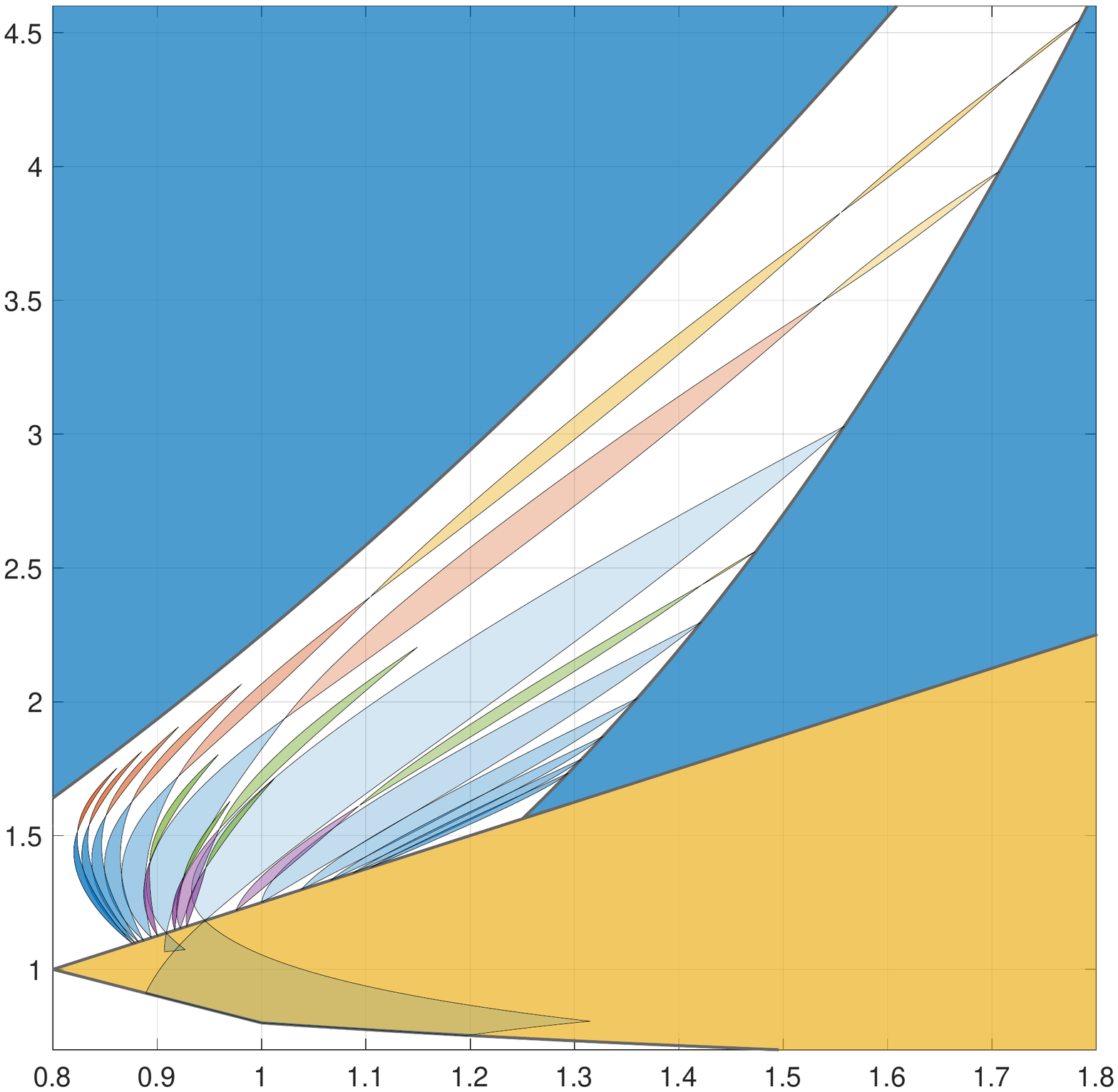}}
\put(120,0){$c_A$}
\put(0,120){$c_B$}

\put(100,20){$A$}
\put(105,75){$BB=D$}
\put(20,90){$AAB=T$}
\put(36,32){\rotatebox{35}{\small $(AAB)(BB)=TD$}}
\put(49,33){\rotatebox{32}{\tiny $TD^2$}}
\put(50,30.6){\rotatebox{30.}{\scalebox{0.5}{\tiny $TD^3$}}}
\put(21,24){\rotatebox{70}{\tiny $T^2 D$}}
\put(17,24){\rotatebox{77}{\scalebox{1.0}{\tiny $T^3 D$}}}
\put(53,57){\rotatebox{36}{\tiny $QTD$}}
\put(33,48.5){\rotatebox{39}{\tiny $QT^2D$}}
\put(63,74){\rotatebox{40}{\tiny $Q^2TD$}}
\put(103,106){\rotatebox{40}{\tiny $Q^3D$}}
\put(115,115){\rotatebox{40}{\tiny $Q^2BD^3$}}

\put(103,97){\rotatebox{40}{\tiny $Q^2D$}}

\put(26,12){\small $TA$}

\put(73,94){\rotatebox{43}{$\nu_4=0$}}
\put(100,40){\rotatebox{18}{$c_Ae_A=c_Be_B$}}
\put(83,50){\rotatebox{54}{$c_A^3e_B=c_Be_A^3$}}

\end{picture}
\end{center}
	\caption{Stability boundaries of various sequences, in $c_A$-$c_B$ parameter space, with $e_A=1$ and $e_B=0.8$. We use the abbreviations $D\equiv BB$, $T\equiv AAB$ and $Q\equiv ABBB$. The blue tongues are regions of the form $(AAB)^{n_1}(BB)^{n_2}=T^{n_1}D^{n_2}$, where from the middle region with $n_1=n_2=1$ (marked $(AAB)(BB)$), $n_1$ increases to the left, and $n_2$  increases to the right. The deepness of the blue colour is proportional to the sequence length. As the parameter $c_B$ is increased through a pinch in each sausage string, an $A$ in the sequence is replaced by a $BB$. A zoom showing more details of a portion of this figure is given in figure~\ref{fig:stab_subcycs_zoom}.
	\label{fig:stab_subcycs}}
\end{figure}

It is clear from the second condition in each list that it is impossible for both $A$ and $B$ to be stable for a single set of parameters. We show the regions of stability of these sequences in the $c_A$-$c_B$ parameter plane in figure~\ref{fig:stab_subcycs}.

\subsubsection{Stability of three-cycles}
\label{sec:stabsigmaT}

As discussed previously, the cycle $\Sigma_T$ has the root sequence $AAB$. To compute the stability of this cycle we need to compute the three matrix products
\[\fl
\mathcal{M}_1=M_{B\rightarrow A} M_{A\rightarrow B}  M_{A\rightarrow A}, \quad \quad  \mathcal{M}_2=M_{A\rightarrow A} M_{B\rightarrow A} M_{A\rightarrow B},  
\quad \quad \mathcal{M}_3=M_{A\rightarrow B}  M_{A\rightarrow A} M_{B\rightarrow A} 
\]
and check when they satisfy the conditions of definition~\ref{lem:stab}. It can be shown then that $\Sigma_T$ is fragmentarily asymptotically stable if the following conditions are satisfied:
\[
\fl
\delta_T\equiv \frac{c_A^2c_B}{e_A^2e_B}>1, \quad \nu_4\equiv-c_B+\frac{(c_A)^2}{e_A}+\frac{c_Ae_A}{e_B}<0, \quad \nu_5\equiv-c_B-\frac{(c_A)^2}{e_A}+\frac{c_Ac_Be_B}{(e_A)^2}<0.
\]
The same conditions can also be computed using different methods, described in~\cite{KS94} and~\cite{Postlethwaite2005}. When the first of these conditions is broken, i.e. $\delta_T<1$, a resonance bifurcation occurs, which, as described in section~\ref{sec:P123}, is degenerate, and does not result in the appearance of a long-period periodic orbit (although we might expect it to with the addition of appropriate fifth order terms to equations~\eref{eq:odes}). 

The remaining two conditions are of a similar type, and correspond to one of the matrices $\mathcal{M}_j$ having an eigenvector with a zero component. In~\cite{KS94,Postlethwaite2005}, the same conditions are arrived at (for very similar systems) by considering perturbations to the cycle $\Sigma_T$  in the $x_4$ and $x_5$ directions respectively. For the parameters we consider in figure~\ref{fig:stab_subcycs}, $\nu_5>0$. The boundary where $\nu_4=0$ is shown. Along this curve, the matrix  $\mathcal{M}_1$ has an eigenvector with a zero in the third component. This means that there are no initial conditions for which a trajectory which performs one $A$ transition, then a $B$ transition, and then an $A$ transition. Or, in terms of the original coordinates, and the cycle $\Sigma_T$, a trajectory which starts on $H^{\iin}_2$ cannot visit $\xi_3$ and $\xi_1$ and then return to $\xi_2$: what is observed is that after the visit to $\xi_1$ the trajectory performs a $B$ transition to $\xi_4$. In the absence of the $x_5$ coordinate (i.e.~in the subspace $P_{1234}$), we would see an exchange of stability between the sequences $AAB$ and $ABBB$.

\section{Numerical results and other cycling patterns}
\label{sec:num}

Although it is in principal possible to generate analytic results about the stability of root sequences of longer length, it is not particularly enlightening (and the expressions are cumbersome at best), and so we proceed to compute the regions of stability for other root sequences numerically. Specifically, for any parameter set, and root sequence $Z$ of length $m$, we can compute the collection of matrices $\mathcal{M}_{Z^j}$ and check the conditions in definition~\ref{defn:fas}. In this section we explain how to extend this notion to numerically determine the stability boundaries in parameter space for the sequence $Z$. 

Definition~\ref{defn:fas} gives a set of conditions which, if satisfied, means that a matrix is of fragmentary asymptotic stability type (abbreviated to simply `stable' in what follows). This results in a dichotomy: either the matrix is stable or it is not. In order to compute stability boundaries using numerical continuation (we use the software MatCont~\cite{matcont}), we require a \emph{continuous} variable which determines the stability, and specifically we define a real-valued scalar,  which is continuous across the stability boundary, and changes sign along the stability boundary. 

From definition~\ref{defn:fas}, we see that a matrix is stable if the eigenvalue with largest absolute value ($\lmax$) lies on the half-line in the complex plane extending from $1$ to $\infty$ along the real axis (which we refer to as $L$), and if all the components of the corresponding eigenvector ($\wmax$) are of the same sign. The matrix can change from stable to unstable in three different ways: (i) by $\lmax$ moving off the half-line $L$; (ii) by one of the components of $\wmax$ changing sign; or (iii) an eigenvalue which does not satisfy the stability conditions becoming larger in magnitude than $\lmax$.

In order to incorporate all these ways in which stability can change, for a matrix $\mathcal{M}$, we define two quantities, $s_{\lambda}(\mathcal{M})$ in terms of the eigenvalues, and $s_w(\mathcal{M})$ in terms of the eigenvectors. We define the set $V$ to be all vectors which satisfy condition~(iii) of definition~\ref{defn:fas}, namely, that all components of the vector are non-zero and of the same sign.

For a matrix $\mathcal{M}$, let the eigenvalue with largest absolute value be $\lmax$. If $\lmax\in\R$,  we label the eigenvalue with second largest absolute value as $\ltwo$, otherwise (that is, if $\lmax$ is one of a complex conjugate pair of eigenvalues), then we label the eigenvalue with third largest absolute value as $\ltwo$. Further, let the eigenvector of $\lmax$ be $\wmax$, with components $\wmax_q$.

We first define, for any eigenvalue $\lambda$ of $\mathcal{M}$, with eigenvector $v_{\lambda}$, the quantity
\[
s_d(\lambda)=\left\{ 
\begin{array}{ll}
|\Im(\lambda)| &  \mbox{if } \Re(\lambda)\geq 1,  \lambda\notin L  \\
|\lambda-1| &  \mbox{if } \Re(\lambda)< 1,  \\
0 & \lambda\in L, v_{\lambda} \in V \\
1 & \lambda\in L, v_{\lambda} \notin V 
\end{array}
\right. 
\]
The first three lines of the definition of $s_d$ give the distance of $\lambda$ from the half-line $L$. The fourth line says that if $\lambda\in L$ but its eigenvector is not in $V$, then $s_d$ is set equal to $1$. We consider those eigenvalues which have $s_d>0$ to be ``unstable'' eigenvalues, and those with $s_d=0$ to be ``stable'' eigenvalues. 

 We then define 
\begin{equation}
s_\lambda(\mathcal{M})=
\left\{
\begin{array}{ll}
|\lmax| - |\lambda_2| & s_d(\lmax)=0 \\
-\min( s_d(\lmax),  |\lmax| - |\lambda_2| ) & s_d(\lmax)>0, s_d(\lambda_2)=0 \\
\green{-s_d(\lmax)} & s_d(\lmax)>0, s_d(\lambda_2)>0 
\end{array}
\right.
\end{equation}
which measures how ``stable'' or ``unstable'' the largest eigenvalue of $\mathcal{M}$ is, giving a positive value when $\lmax$ is a ``stable'' eigenvalue, and a negative value when $\lmax$ is ``unstable''.
Specifically, we have, in the first line, the difference in absolute values of the largest and second largest eigenvalues, in the case where $\lmax$ is stable. If $\lmax$ is unstable but $\lambda_2$ is stable, then the second line gives (minus) the smaller of $s_d(\lmax)$ and the difference in absolute values of the largest and second largest eigenvalues. When both $\lmax$ and $\lambda_2$ are unstable, the third line gives (minus) the smaller of the function $s_d$ evaluated at both eigenvalues.

We then define
\[
s_w(\lambda)=\left\{ 
\begin{array}{ll}
\min_{q,l}(({v_\lambda})_q({v_\lambda})_l)  &  \mbox{if } \lambda\in L  \\
-1 &  \mbox{if } \lambda \notin L,  
\end{array}
\right. 
\]
The quantity $s_w$ is only positive when all components of the eigenvector have the same sign, and $\lambda\in L$. It is zero if one of the components of the eigenvector is zero (and $\lambda\in L$). Otherwise it is negative. It will change sign in a continuous fashion for $\lambda\in L$ when one of the components of the eigenvector changes sign. Similarly to $s_{\lambda}$,  those eigenvalues which have $s_w>0$ are ``stable'' eigenvalues, and those with $s_w<0$ are ``unstable'' eigenvalues.

Finally, we define
\begin{equation}\label{eq:s}
s(\mathcal{M})=\left\{
\begin{array}{ll} 
s_\lambda(\mathcal{M})s_w(\lmax) & s_\lambda(\mathcal{M})>0, s_w(\lmax)>0 \\ 
-|s_\lambda(\mathcal{M})s_w(\lmax)| & \mathrm{otherwise.}
\end{array}
\right.
\end{equation}

This gives a positive number when both $s_{\lambda}$ and $s_w$ are positive, and a negative one otherwise, and is continuous across the boundary where $s=0$.
For a collection of matrices $\{\mathcal{M}_{Z^j}, j=1,\dots m\}$ we compute $s(\mathcal{M})$ for all matrices in the collection, and take the minimum, in order to compute the stability of the root sequence $Z$.

 Note that it is possible for $s(\mathcal{M})$ to be equal to zero when there is not a change in stability (specifically, in the case where $|\lmax|=|\ltwo|$ and  both $\lmax$ and $\ltwo$ lie on $L$), but it will not change sign.

\begin{figure}
\begin{center}
\setlength{\unitlength}{1mm}
\begin{picture}(130,130)(0,0)
\put(0,0){\includegraphics[trim= 0.cm 5cm 0cm 4.5cm,clip=true,width=130mm]{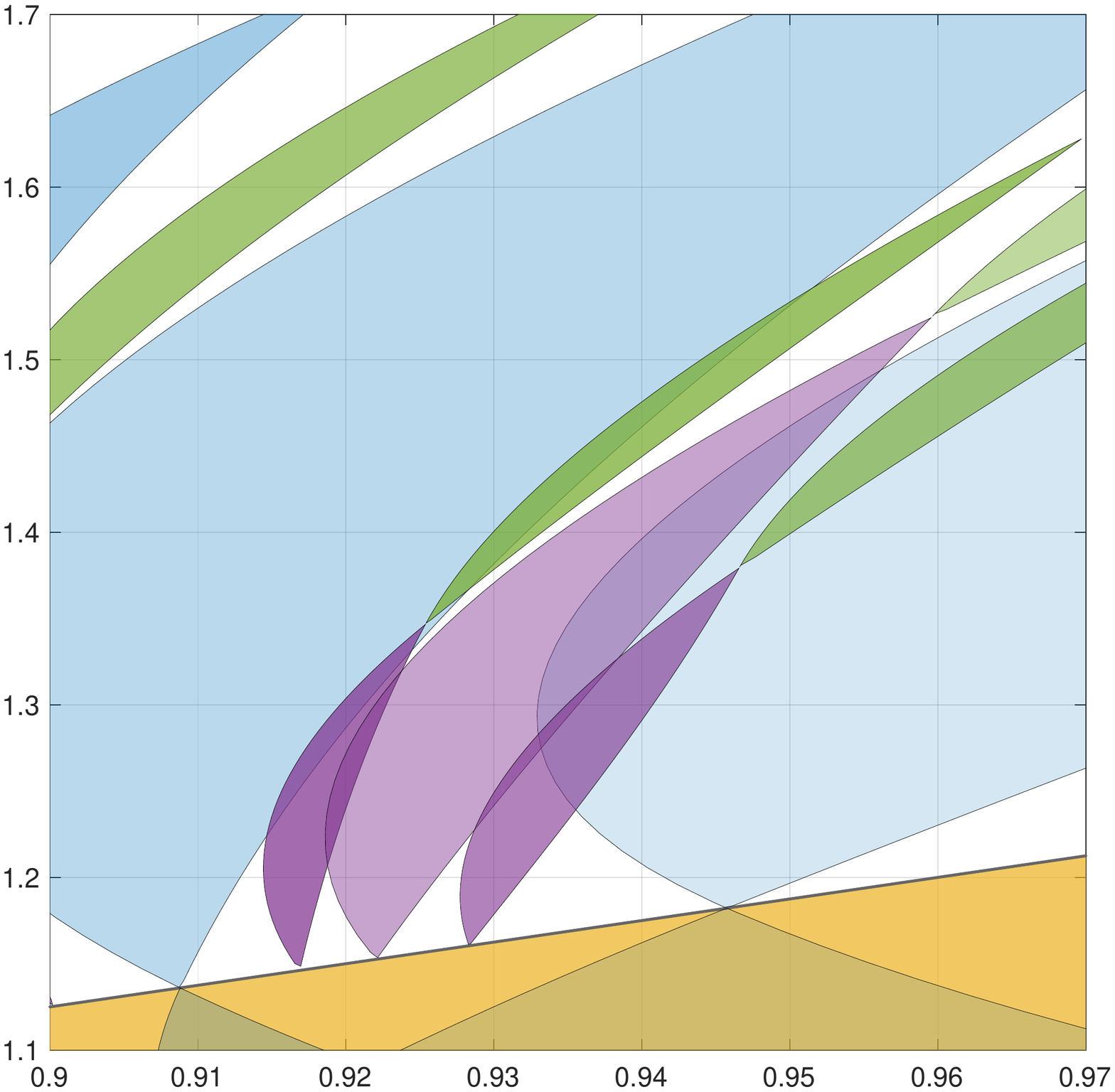}}
\put(120,0){$c_A$}
\put(0,120){$c_B$}

\put(115,17){$A$}
\put(15,40){\rotatebox{60}{ $T^2D$}}
\put(50,40){\rotatebox{45}{ $T^2DTD$}}
\put(90,40){\rotatebox{30}{ $TD$}}

\put(34,20){\rotatebox{75}{\tiny $T^2DT^2DTD$}}
\put(83,80){\rotatebox{37}{\tiny $QTDT^2DTD$}}
\put(116,93.2){\rotatebox{35}{\tiny $QTDTD$}}

\put(57,23){\rotatebox{52}{\tiny $T^2DTDTD$}}
\put(100,71){\rotatebox{33}{\tiny $QTDTDTD$}}

\put(83,10){$TA$}
\put(24,8){\tiny $T^2A$}



\end{picture}
\end{center}
	\caption{The figure shows a zoom of a portion of figure~\ref{fig:stab_subcycs}. Recall the abbreviations $D\equiv BB$, $T\equiv AAB$ and $Q\equiv ABBB$. 
		\label{fig:stab_subcycs_zoom}}
\end{figure}

Using this algorithm, we are able to find a series of `tongues' of stability regions of different sequences, arranged into `strings of sausages', thirty-four of which are shown in figures~\ref{fig:stab_subcycs} and~\ref{fig:stab_subcycs_zoom}. These sausages are reminiscent of the shape of resonance tongues in piecewise smooth systems~\cite{wei1987,campbell1996,szalai2009,simpson2016border,simpson2018structure}. This is not overly surprising, because the map which we use to derive these stability boundaries, is of course, defined in a piecewise fashion. However, in those systems, in each string of sausages, the same object is stable, but this is not the case here.

%
%
%
%
%
%
%

We observe that the strings of sausages have the following properties:
\begin{itemize}
\item All but two of the strings occur entirely between the curves $c_Ae_A=c_Be_B$ (the upper stability boundary of the $A$ sequence), the curve $c_A^3e_B=c_Be_A^3$ (the left-most stability boundary of the $B$ sequence), and the curve $\nu_4=0$ (the right-hand stability boundary of the sequence $AAB$).
\item The vast majority of the stable sequences are made up only of the components $(BB)$, $(AAB)$ and $(ABBB)$. For clarity, we relabel these sequences as $D$, $T$ and $Q$ respectively.
\item For tongues which abutt the curve $c_Ae_A=c_Be_B$, between any two tongues with sequences $S_1$ and $S_2$, there is another tongue which has the sequence $S_1S_2$, giving a Farey-like sequence. In particular, we predict the existence of stable regions for all sequences of the form $(AAB)^{n_1}(BB)^{n_2}$, $n_1,n_2\in\N$. We found 12 of these numerically, and they are coloured blue in figure~\ref{fig:stab_subcycs}.
Between each of these blue-coloured tongues, there are additional tongues --- see the zoomed figure~\ref{fig:stab_subcycs_zoom} --- here coloured purple. Notice the large purple tongue labelled $T^2DTD$ between the blue regions labelled $T^2D$ and $TD$. Similarly, there are smaller purple tongues which have the sequence of the adjacent tongues.
\item As the parameter $c_B$ is increased through a `pinch' in the sausage string, an $A$ in the sequence transforms into a $BB$. Normally, this results in a $T$ transforming into a $Q$. For instance, notice the string with sequences $T^2D$, $TQD$ and $Q^2D$ in figure~\ref{fig:stab_subcycs}.
\end{itemize}

Figure~\ref{fig:stab_subcycs_zoom} shows a zoom showing a selection of additional tongues abutting the line $c_Ae_A=c_Be_B$. Notice that the middle purple area $T^2DTD$ appears between $T^2D$ and $TD$, and the other two purple regions have sequences formed by a concatenation of their neighbours. In this figure, we can also see two stability regions which occur entirely within the $A$ stability region, on the other side of the line $c_Ae_A=c_Be_B$. These regions have sequences $AAAB$ and $TAAAB$. Each of these touches the line $c_Ae_A=c_Be_B$ at the same point as a region on the other side of that line with a sequence with an $A$ replaced by a $BB$, thus continuing the pattern noted above. We were unable to find any other such regions within the $A$ stability region.

The sausage strings appear to have a number of other interesting properties. For instance, the upper boundaries of the $T^{n_1}D$ tongues sit on a straight line, as do the upper boundaries of the $QT^{n_1}D$ tongues. The number of tongues in each sausage string is always is always odd, and (aside from those of type $TD^n$) increase with the position of the original tongue in the `Farey sequence'.

The region of parameter space shown in figure~\ref{fig:stab_subcycs} is mostly (but not all) within the region of parameter space for which the sufficient conditions for asymptotic stability of the network $\Sigma$ computed in 
~\cite{Podvigina2020} and~\cite{Afraimovich2016} apply (namely, that $\min(c_A, c_B)>\max(e_A, e_B)$; recall the grey shaded region in figure~\ref{fig:stab_bounds}). We conjecture that the remaining white space between the $A$, $AAB$ and $BB$ regions in figure~\ref{fig:stab_subcycs} is actually completely filled with (overlapping) stability tongues.

\subsection{Loss of stability of sequences}

As well as computing the boundaries of stability for each of the tongues shown in figure~\ref{fig:stab_subcycs}, we can further identify how the stability is lost. In all cases we have observed, stability is lost through the breaking of condition (iii) of definition~\ref{defn:fas}, that is, one of the matrices in the collection for that sequence has an eigenvector with a zero component. Furthermore, in all the cases we have checked, it is the third component of the eigenvector which is zero. As we noted at the end of section~\ref{sec:tmatrices}, identifying which of the matrices $\mathcal{M}_{Z^j}$ has the eigenvector with a zero eigenvalue can tell us what happens when stability is lost. We demonstrate this using the stability tongue for the sequence $T^2 D\equiv AABAABBB$ as an example.

On the left hand side of the stability tongue for the sequence $T^2D$, we find that the matrix 
\[
\mathcal{M}_8\equiv M_{A\rightarrow A} M_{B\rightarrow A}  M_{B\rightarrow B}M_{B\rightarrow B} M_{A\rightarrow B}  M_{A\rightarrow A}M_{B\rightarrow A} M_{A\rightarrow B}
\]
has an eigenvector with a zero in the third component. This means that there are no initial conditions which allow trajectories to complete this sequence in the order $BAABBBAA$. In particular, trajectories which perform the first seven of these transitions would then perform a $B$, rather than an $A$, as the next transition in the sequence. If we assume, for now, that the final $A$ is actually replaced in the sequence by a $BB$, then we get the sequence $BAABBBABB$, which reordered is equivalent to $ABBBAABBB\equiv QTD$, which is the sequence observed in the next tongue up the sausage string.

On the right hand side of the stability tongue for the sequence $T^2D$, we find that the matrix 
\[
\mathcal{M}_6\equiv  M_{B\rightarrow A} M_{A\rightarrow B}M_{A\rightarrow A} M_{B\rightarrow A} 
M_{B\rightarrow B}M_{B\rightarrow B} M_{A\rightarrow B}  M_{A\rightarrow A}\]
has an eigenvector with a zero in the third component. This is different to the left hand boundary. Here, it means that there are no initial conditions which allow trajectories to complete the sequence in the order $ABBBAABA$. Taking a similar approach, and replacing the final $A$ in this sequence with a $BB$, gives the sequence $ABBBAABBB\equiv QTD$, again. 

We see a similar pattern arise on the boundaries of the other tongues: namely, that it is an $A$ transition which is the first transition which cannot occur once the boundary has passed, and if this $A$ is replaced in the sequence with a $BB$, then the resulting sequence is the one which is found up one level in the sausage string. A detailed analysis of why this occurs is beyond the scope of this paper, but will be investigated in future work.

\section{Transitions to irregular cycling behaviour}
\label{sec:irreg}

Figure~\ref{fig:tslog}(b) showed an example of a trajectory which approached the network, but 
 in a complicated manner, that is, the visits to the equilibria were a non-repeating sequence. Behaviour of this kind was also noted for a network between six equilibria in~\cite{Postlethwaite2005}, and there has been other work on the existence of such switching behaviour in other types of networks~\cite{aguiar2004dynamics,homburg2010switching,Kirk2010,castro2016switching}.

We identify a region of parameter space where this type of behaviour appears to be typical: we can find no stable regular sequences here.
In figure~\ref{fig:stab_bounds}, this parameter region is the area coloured white in the lower right corner. Note that this region is not within the region of parameter space where the sufficient condition for asymptotic stability of $\Sigma$ from~\cite{Podvigina2020,Afraimovich2016} applies. However, we note from numerical simulations that the network appears to be highly attracting here. We observe two different ways in which this irregular behaviour arises. Firstly, as a consequence of quasi-periodic behaviour (resulting from bifurcations from the equilibrium $\xi_Q$) which gets closer to the network. Secondly, as an instability of the heteroclinic cycle $\Sigma_{TQ}$. We present numerical evidence of both of these mechanisms below, but leave a detailed investigation of this behaviour to later work.

\subsection{Irregular behaviour arising from quasi-periodic solutions}

\begin{figure}
\begin{center}
\setlength{\unitlength}{1mm}
\begin{picture}(130,90)(-5,3)
\put(0,45){\includegraphics[trim= 2.8cm 5cm 2.3cm 6.5cm,clip=true,width=130mm]{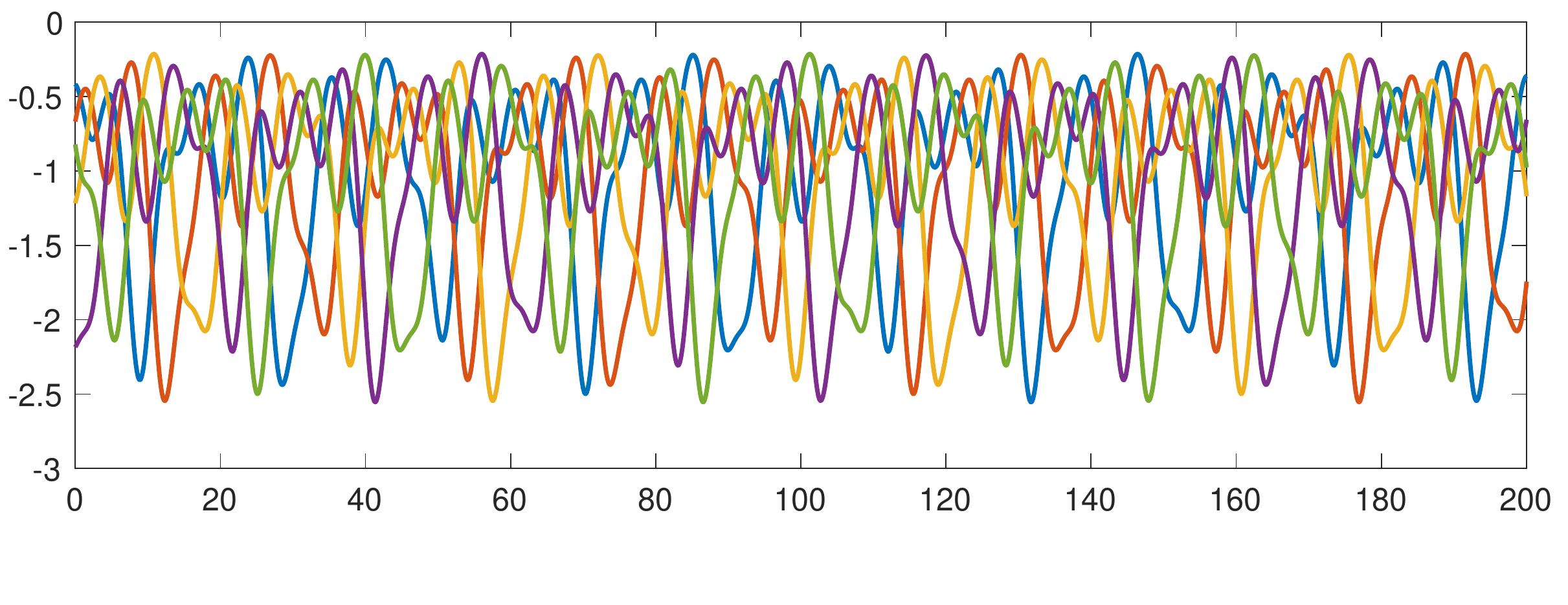}}
\put(0,0){\includegraphics[trim= 2.8cm 5cm 2.3cm 6.5cm,clip=true,width=130mm]{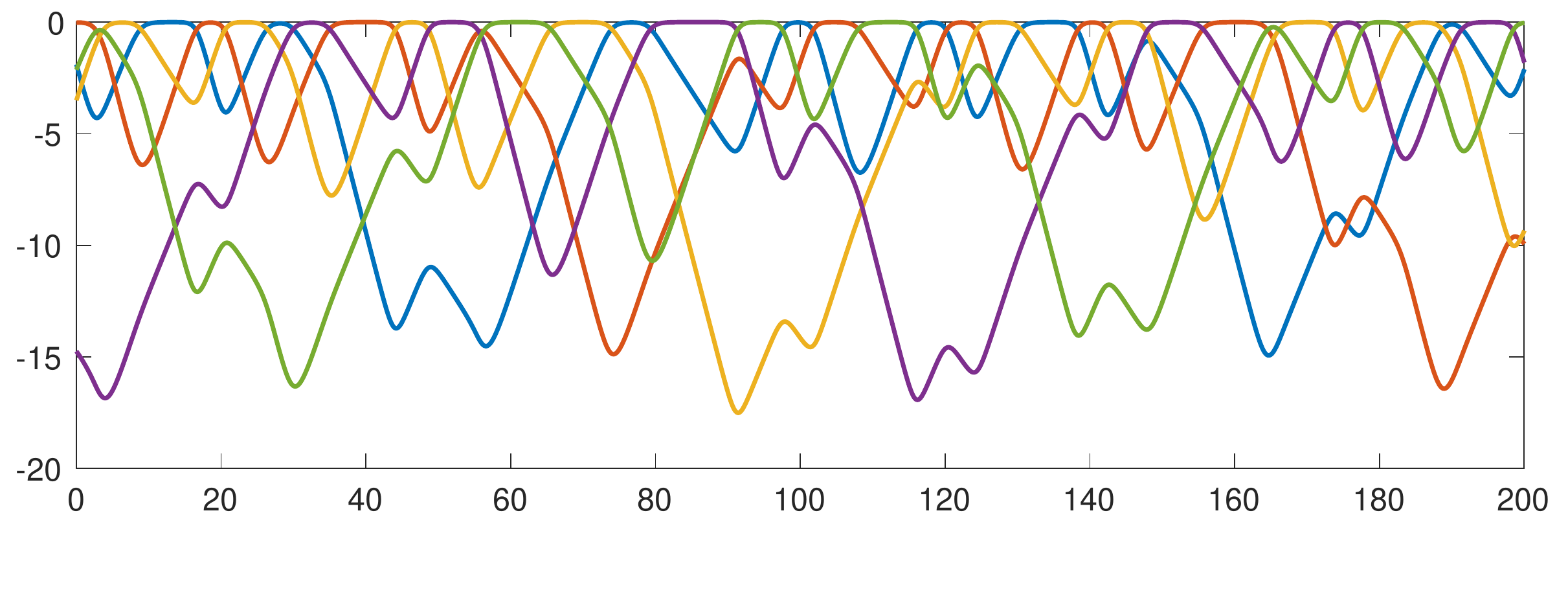}}

\put(-3,95){(a)}
\put(-3,50){(b)}
\put(-5,70){\rotatebox{90}{$\log(x_j)$}}
\put(-5,25){\rotatebox{90}{$\log(x_j)$}}
\put(130,4){$t$}

\end{picture}
\end{center}
\caption{The figures show typical time series of equations~\eref{eq:odes}, on a logarithmic scale. The lines coloured blue, red, yellow, purple and green are the logarithm of the  coordinates $x_1,\dots,x_5$ respectively. Parameters are: panel (a): $c_A=1.2$, $c_B=0.5$; panel (b): $c_A=1.25$, $c_B=0.5$; $e_A=1$ and $e_B=0.8$ throughout. 
	\label{fig:tsquasi}}
\end{figure}
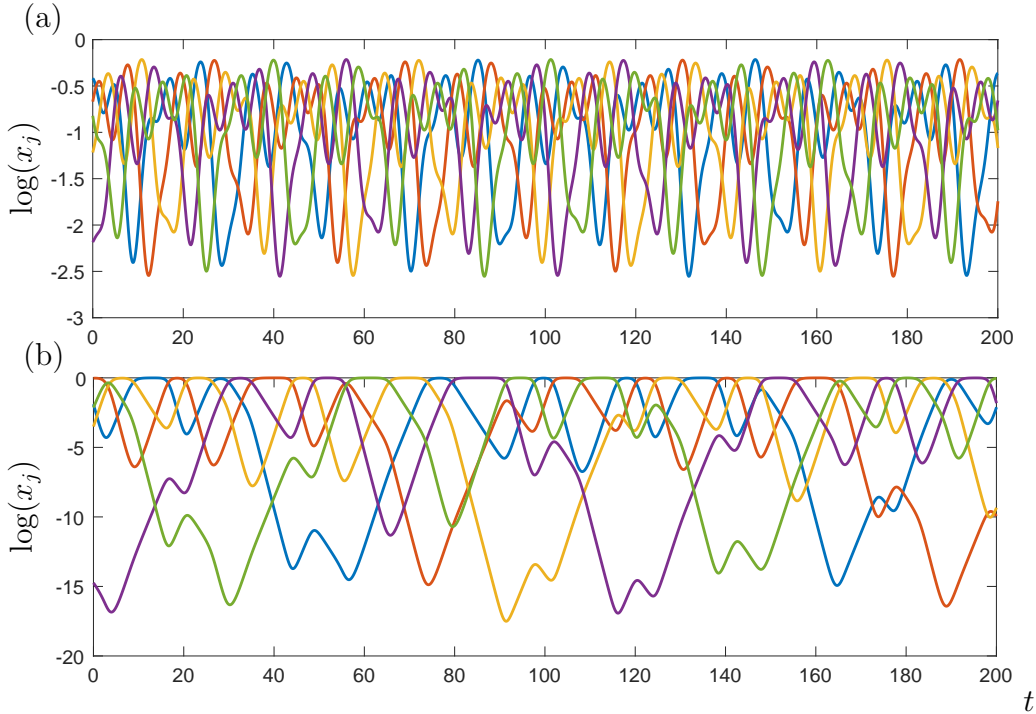

In figure~\ref{fig:ts}(d) we noted that following the Hopf bifurcation from $\xi_Q$, we observe quasi-periodic behaviour at $c_A=1.02$, $c_B=0.5$. In figure~\ref{fig:tsquasi} we show two further timeseries for increasing value of $c_A$, both in logarithmic coordinates. As $c_A$ is increased, the quasiperiodic behaviour grows in amplitude, but due to the existence of the invariant sphere, it becomes essentially `trapped' between the equilibria of the network $\Sigma$. As a result, the trajectory spends increasingly long times in neighbourhoods of the equilibria. In figure~\ref{fig:tsquasi}(b) we can see clear visits of the trajectory to the five equilibria $\xi_j$, during which all but one of the coordinates are very small. However, the value of the smallest of these coordinates, which can be thought of as an approximate distance of the trajectory from $\Sigma$, does not tend to decrease over time, as can be observed in both trajectories in figure~\ref{fig:tslog}.
As $c_A$ is increased further, we observe that the distance of the irregular trajectories from the network decreases, although the network does not appear to be `attracting' until $c_A$ is greater than about $1.3$, when we see a sustained approach of the trajectory towards the network.

\subsection{Irregular behaviour arising from a bifurcation of $\Sigma_{TQ}$}

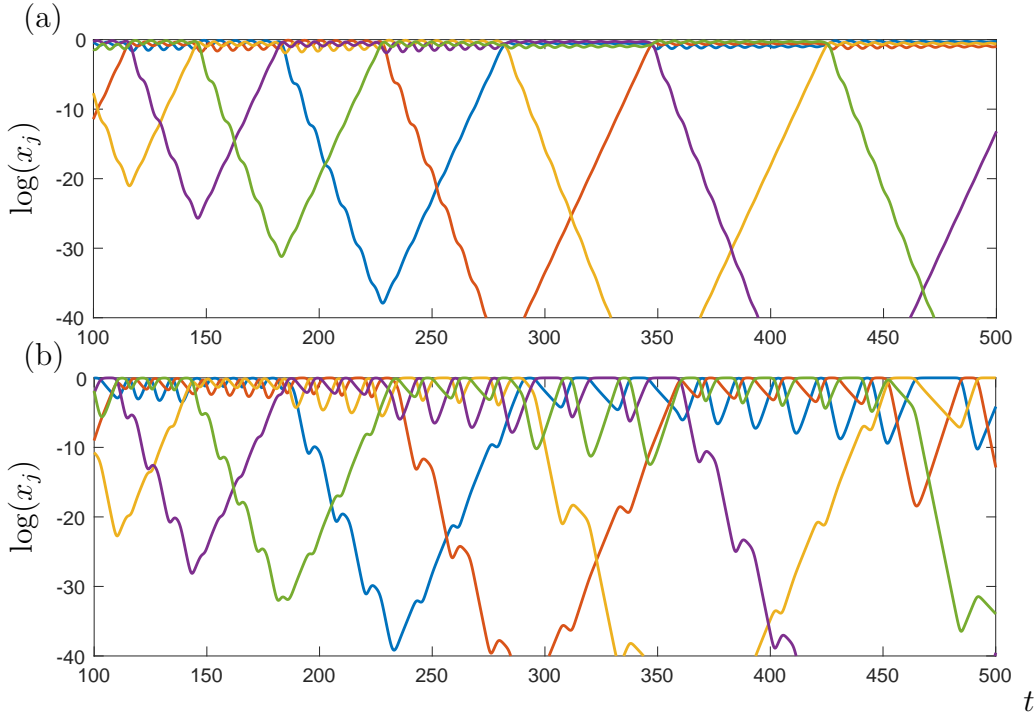
\begin{figure}
\begin{center}
\setlength{\unitlength}{1mm}
\begin{picture}(130,90)(-5,3)
\put(0,45){\includegraphics[trim= 2.8cm 5cm 2.3cm 6.5cm,clip=true,width=130mm]{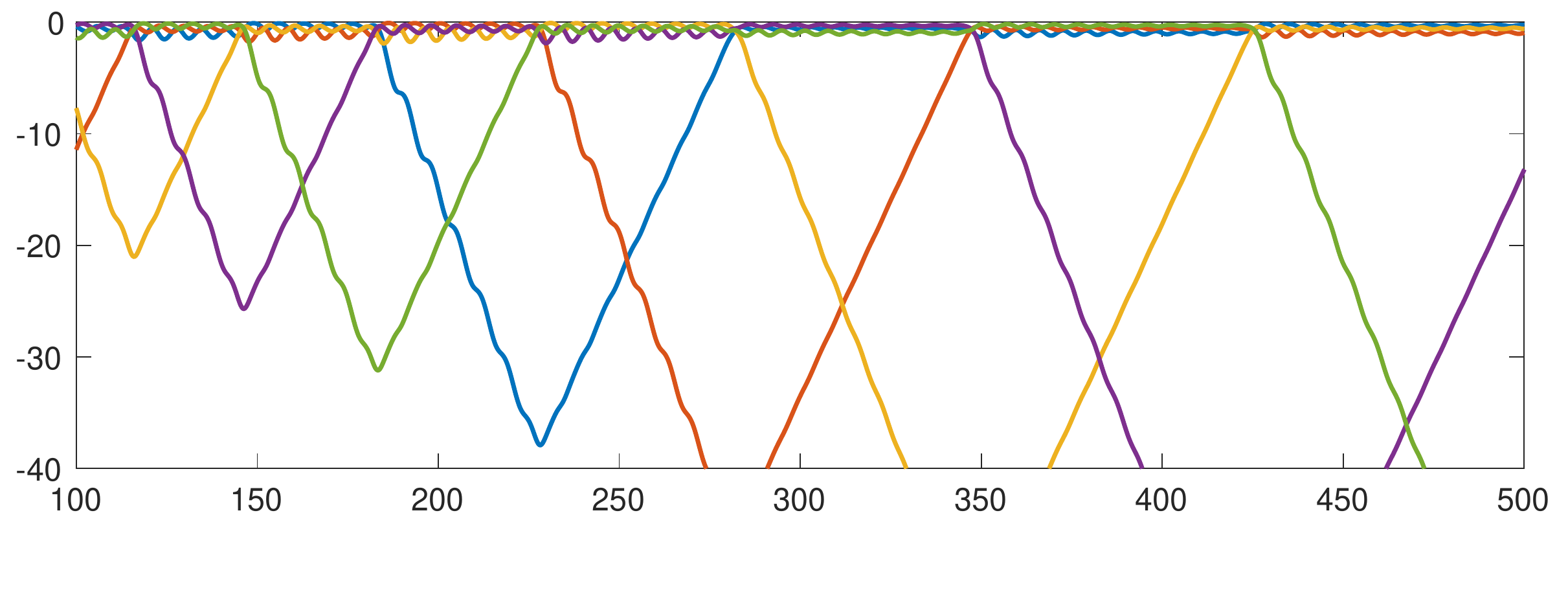}}
\put(0,0){\includegraphics[trim= 2.8cm 5cm 2.3cm 6.5cm,clip=true,width=130mm]{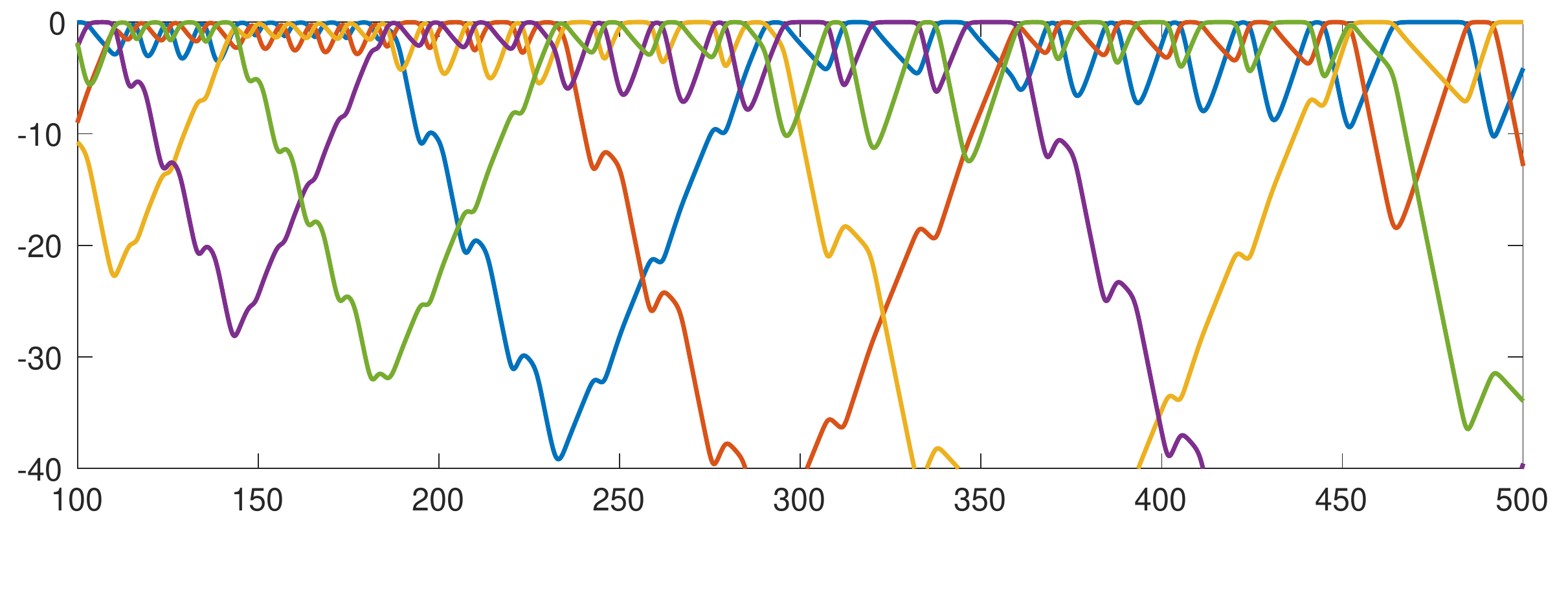}}

\put(-3,95){(a)}
\put(-3,50){(b)}
\put(-5,70){\rotatebox{90}{$\log(x_j)$}}
\put(-5,25){\rotatebox{90}{$\log(x_j)$}}
\put(130,4){$t$}

\end{picture}
\end{center}
\caption{The figures show typical time series of equations~\eref{eq:odes}, on a logarithmic scale. The lines coloured blue, red, yellow, purple and green are the logarithm of the  coordinates $x_1,\dots,x_5$ respectively. Parameters are: panel (a): $c_A=1.6$, $c_B=0.25$, and so $\delta_T>1$; panel (b): $c_A=1.65$, $c_B=0.35$, and so $\delta_T<1$; $e_A=1$ and $e_B=0.8$ throughout. 
	\label{fig:tssigmatq}}
\end{figure}

In figure~\ref{fig:ts}(b) we show a trajectory approaching the heteroclinic cycle $\Sigma_{TQ}$ (defined in section~\ref{sec:P123}). Rather than approaching the equilibria on the coordinate axes, this cycle approaches equilibria which have three non-zero coordinates, and lie in the interior of the subspaces $P_{j,j+1,j+2}$. For the dynamics restricted to these three-dimensional subspaces, these equilibria loose stability when $\delta_{T}$ is increased through $1$, shown by a blue curve in figure~\ref{fig:stab_bounds}. In figure~\ref{fig:tssigmatq} we show timeseries for two parameter sets, on each side of this line, with $c_A=1.6$. Specifically, in panel (a), we can see that within each visit to a subspace  $P_{j,j+1,j+2}$ (i.e. where two of the coordinates are very small), the oscillations are decaying: the trajectory is approaching the equilibrium in the interior. In panel (b), we see the opposite: during each visit to a subspace  $P_{j,j+1,j+2}$, the amplitude of the oscillations increase, and in fact we see a sequence of visits to the three equilibria $\xi_j$, $\xi_{j+1}$ and $\xi_{j+2}$. However, the number of visits made to these equilibria before we switch to the next subspace is not the same each time. This behaviour was termed \emph{irregular cycling} in~\cite{Postlethwaite2006}, and in that paper we proved its existence for an open region of parameter space. In this system, this particular type of irregular behaviour seems to be restricted to a fairly small region of parameter space: it rapidly breaks down into the much more irregular behaviour of the type seen in figure~\ref{fig:tslog}(b).

\section{Discussion}
\label{sec:disc}

In this paper we have developed a method for determining regions of parameter space in which different patterns of approaching a heteroclinic network can be found, and have applied this to a model of five-species cyclic competition. We find some complicated and intriguing patterns of stability regions in parameter space, which are reminiscent of those found in other piecewise studies. The numerical technique we have developed here could easily be used to analyse the stability of sequences of visited equilibria in other heteroclinic networks with two-dimensional unstable manifolds.

This work offers several potential avenues for further study. Firstly, we would like to be able to prove some of the observations made about the patterns of stability tongues found in section~\ref{sec:num}. This would likely involve further analysis of the \Poincare maps~$\Phi_B$~\eref{eq:phib} and~$\Phi_A$~\eref{eq:phia}, before the approximations are made which assume we are close to the one-dimensional connections in $\hat{\Sigma}$. From the numerical results, it appears that the transitions between the tongues in each string of sausages is associated with a change from a $A$ pattern to an $BB$ pattern, which is exactly when the above approximation is not valid.
Secondly, our calculations tell us only the regions where sequences have a basin of attraction of positive measure; we could further compute the shape of the $\delta$-local basin of attraction (see equation~\eref{eq:deltabasin}) using the `stability index'~\cite{Podvigina2011,garrido2019stability}. We have done these calculations \green{numerically}, and found that except for the $A$ sequences, all sequences have a `cusp' shaped local basin - that is, one in which the measure decreases to zero as the network is approached. \green{These numerical calculations are supported by analytical calculations in~\cite{castro2021stability}, for the $A$, $B$ and $AAB$ cycle}. However, these sequences are not hard to find by randomly selecting initial conditions for numerical integrations, which seems somewhat counter-intuitive. In fact, what we observe is that many trajectories wander away from a particular sequence for some transient period before settling down, meaning that the entire basin of attraction is much larger than predicted using only a local analysis. This type of behaviour was noted in both~\cite{KS94} and~\cite{Postlethwaite2005}; in the latter the behaviour was termed \emph{essential quasi-asympototic stability}.
Thirdly, the numerical simulations clearly indicate that the network has strong attractive properties in regions of parameter space where there are no stable root sequences, and where the sufficient conditions for stability derived in~\cite{Podvigina2020,Afraimovich2016} do not apply, namely the region where we observe irregular cycling. It would be of great interest to investigate this behaviour further. 

It is also natural to ask the question: are the sorts of dynamics we observe for the Rock-Paper-Scissors-Lizard-Spock network typical for larger networks, or for networks in which the $\Z_5$-symmetry is broken? To first address the issue of broken symmetry (i.e.~there is a different set of eigenvalues at each equilibria, albeit with the same signs): we believe that much of the observed dynamics would remain similar. Numerical simulations (by integration) indicate that complicated dynamics are still possible. Although we have not computed any stability boundaries for the broken symmetry case, the technique would be exactly the same, only with more book-keeping: one would now need to keep track of equilibria visited in addition to whether $A$ or $B$ type connections were traversed.

We also expect to see similarly complicated dynamics in larger networks, with the added complication that there would now be more than one possible network topology, even if we consider only networks between equilibria with a single species present, and retain the restriction 
\green{that there is a symmetry between the equilibria which preserves the network structure (that is, all the equilibria are the `same'). It is possible to use group-theoretic methods to compute all the possible network topologies for `small' number of equilibria $k$ (as in, e.g.~\cite{potovcnik2014groups,potovcnik2014census,holt2020census}), and the number grows very rapidly with $k$.}
As for the non-symmetric case, the technique we have developed in this paper for analysing the behaviour could be applied to any of these networks. A full cataloguing of the behaviour for these networks will involve classifying the networks by some measures of their topology. Initial investigations we have made have shown that networks with $k$ odd, which contain a $\Delta$-clique (as the Rock-Paper-Scissors-Lizard-Spock network does), have a very similar pattern of tongues of stability regions arranged into strings of sausages. Understanding why this behaviour occurs would generalise the $A\rightarrow BB$ transition discussed above.

\section*{Acknowledgements}

This project began during less traumatic times when international travel was still possible, and we acknowledge the London Mathematical Society for financial support through a Research in Pairs (Scheme 4) grant, and the hospitality of and financial support from the Departments of Mathematics at both the University of Auckland and the University of Leeds. The majority of this research was done during COVID lockdown in 2020, and both authors are grateful to their partners and children for giving them time and space to work. We are grateful to Bernd Krauskopf for helpful conversations about the stability `sausages', to Gabriel Verret for \green{discussions on how to compute the possible symmetric graphs with larger numbers of equilibria mentioned in the discussion}, and for some constructive comments from anonymous referees.
CMP is grateful for additional support from the Marsden Fund Council from New Zealand Government funding, managed by The Royal Society Te Ap\={a}rangi, and from the London Mathematical Laboratory.

\section*{References}

\bibliographystyle{unsrt}
\bibliography{clairebib}

\end{document}